\documentclass{birkmult}
 \newtheorem{thm}{Theorem}[section]

 \newtheorem{prop}[thm]{Proposition}
 \theoremstyle{definition}
 \newtheorem{defn}[thm]{Definition}
 \theoremstyle{remark}
 \newtheorem{rem}[thm]{Remark}
 
 \numberwithin{equation}{section}

 \usepackage{amsmath,amsthm}
\usepackage{amsfonts,amssymb}
\usepackage{xypic}

\newcommand{\uu}{\underline}

\newcommand{\ki}{\Omega^ {\mathcal K}}
\newcommand{\hi}{\Omega^ {\mathcal H}}
\newcommand{\oi}{\Omega^ {\mathcal P}}

\newcommand{\m}{\mathcal}

\begin{document}
%
%
%
%
%
%
%
%
%
\title[Generalized Repeated Interaction Model and Transfer Functions]
 {Generalized Repeated Interaction Model and Transfer Functions}
\author{Santanu Dey}

\address{
Department of Mathematics,\\
Indian Institute of Technology Bombay,\\
Powai, Mumbai- 400076, India}

\email{dey@math.iitb.ac.in}

\author{Kalpesh J. Haria}
\address{Department of Mathematics,\\
Indian Institute of Technology Bombay,\\
Powai, Mumbai- 400076, India}
\email{kalpesh@math.iitb.ac.in}


\date{July, 2013}

\begin{abstract}
Using a scheme involving a lifting of a row contraction we
introduce a toy model of repeated interactions between quantum
systems. In this model there is an
outgoing Cuntz scattering system involving two wandering
subspaces. We associate to this model an input/output linear system which leads
to a transfer function. This transfer function is a multi-analytic
operator, and we show that it is inner if we assume that the
system is observable. Finally it is established that transfer functions coincide
with characteristic functions of associated liftings.
\end{abstract}

\maketitle

\section{\textbf Introduction}

In page  287 of the article \cite{Go10} the author has commented the following while comparing \cite{Go10} with
\cite{DG07,DG11}: In \cite{DG07} a row contraction $\uu A$ on a Hilbert space $\m H$ with a one-dimensional eigenspace
is considered and the theory of minimal isometric dilations is used. The characteristic function introduced in \cite{DG11} is
a multi-analytic operator associated to a lifting and the ergodic case is studied in detail in \cite{DG07}. In \cite{Go10} minimality is not considered but one starts
with an interaction $U$ (which is a unitary operator) in a scheme similar to \cite{DG07} and and obtains a multi-analytic operator which represents the transfer function of an input-output
system associated with the interaction. It is expected that the scheme developed \cite{Go10} is more directly applicable to physical models.
In the setting of \cite{DG11} the assumption of a one-dimensional eigenspace is dropped
and the theory is much more general in another direction. A further integration of these schemes in the future may help to
remove unnecessarily restrictive assumptions
of the toy model considered in \cite{Go10} and lead to the study of other and of more realistic models.

 This paper achieves some of these objectives. In the
model of repeated interactions between quantum systems, also called a noncommutative Markov chain, studied in
\cite {Go10} (cf. \cite {Go04}) for given three Hilbert spaces $\m
H, \m K$ and $\m P$ with unit vectors $\hi, \ki$ and $\oi$ an {\em
interaction} is defined to be a unitary operator $U : \m H \otimes
\m K \to \m H \otimes \m P$ such that
\begin{equation}
U(\hi \otimes \ki) = \hi \otimes  \oi.
\end{equation}
Define $\m K_\infty := \bigotimes_{i=1}^\infty  \m K$ and $ \m
P_\infty := \bigotimes_{i=1}^\infty  \m P$ as infinite tensor
products of Hilbert spaces with distinguished  unit vectors. We denote $m$-th
copy of $\m K$ in $\m K_\infty$ by $\m K_m$ and set $\m K_{[m,n]} :=
\m K_m \otimes \cdots \otimes \m K_n$ for $m \leq n$. Similar notations
are also used with respect to $\m P$.
The repeated interaction is defined as
\[U(n):=U_n \ldots U_1 : \m H \otimes \m K_\infty \to \m H \otimes \m P_{[1,n]} \otimes \m
K_{[n+1,\infty)}\] where $U_i$'s are copies of $U$ on
the factors $\m H \otimes \m K_i$ of the infinite tensor products and $U_i$'s  leaves other factors fixed.
Equation (1.1) tells us that the tensor product of the vacuum vectors $\hi, \ki$ (along with $\oi$)  represents a state of the coupled system which is not affected by the
interaction $U.$ This entire setting represents interactions of an
atom with light  beams  or fields. In particular $\hi$ in
\cite{Go10} is thought of as the vacuum state of an atom, and $\ki$
and $\oi$ as a state indicating the absence of photons.

In the generalized repeated interaction model that we introduce in this
article we use a pair of  unitaries to encode the interactions
instead of one  unitary as follows:\\
Let $\tilde{\m H}$ be a (closed) subspace of $\m H,$ and $U : \m H\otimes\m K \to \m H \otimes \m P$ and  $\tilde U :
\tilde{\m H} \otimes\m K \to \tilde{\m H} \otimes \m P$ be two
unitaries such that
\begin{equation}
 U(\tilde h \otimes \ki) = \tilde U(\tilde h \otimes \ki)~\mbox{for all }~\tilde h \in \tilde{\m H}.
\end{equation}
We fix $\{\epsilon_1,\ldots, \epsilon_d \}$ to be an orthonormal basis of
$\m P.$ The equation (1.2) is the analog of the equation (1.1) for
our model and thus our model can be used for the setting where a quantum system
 interacts with a stream of copies of another quantum system in such a way that there is no backaction (so we get a Markovian type of dynamics)
and such that there is a certain kind of subprocess. In the model of \cite{Go10} the vacuum state $\hi$ of an atom plays an important role.
For a model describing interaction of a quantum system with a stream of copies of another quantum sytem
we need that the computations do not involve any fixed unit vector $\hi$ and we are able to achieve this  in our model by using a pair of
unitaries. Instead of $\hi$ we now have a kind of subprocess, described by $\tilde {U},$ which can be treated on the same level as the full process, described by $U.$

The main condition imposed on the unitary $U: \m H\otimes\m K \to \m H \otimes \m P$ in order to get a generalized
interaction model is that
$U (\tilde{\m H} \otimes \ki) \subset \tilde{\m H} \otimes \m P$
(cf. Proposition 3.1 of \cite{Go11} for an interesting consequence of this assumption).
We can then define $\tilde{U}$ restricted to $\tilde{\m H} \otimes \ki$ as $U$ restricted to $\tilde
{\m H} \otimes \ki,$ and assume that $\m H \otimes \m P$ is big enough to
allow a unitary extension $\tilde{U}: \tilde{\m H} \otimes \m K \rightarrow \tilde{\m H} \otimes \m P.$
The focus of the study done here, as also in \cite{Go10}, is to
bring out that certain multi-analytic operators of the multivariate
operator theory are associated to noncommutative Markov chains and
related models, and that these operators can be exploited as powerful
tools. These operators occur as central objects in various context
such as in the systems theory related works (cf. \cite{BV05}) and
noncommutative multivariable operator theory related works
(cf. \cite{Po89b}, \cite{Po95}).

A tuple $\underline T= (T_1,\ldots,T_d)$ of operators $T_i$'s on a
common Hilbert space $\m L$ is called a \textit{row contraction}
if $\sum_{i=1}^{d}T_iT_i^* \leq I.$ In particular if
$\sum_{i=1}^{d}T_iT_i^* = I$, then the tuple $\underline T = (T_1,
\ldots, T_d)$ is called \textit{coisometric}. We introduce the notation $\tilde {\Lambda}$ for the free semigroup with generators $1,\ldots, d$. Suppose $T_1, \ldots ,T_d \in \m B(\m L)$ for a Hilbert space $\m L.$
If $\alpha \in \tilde \Lambda $ is the word $\alpha_1\ldots\alpha_n$ with length  $|\alpha |= n$, where each $\alpha_j \in \{1,\ldots,d\},$
then $T_\alpha$  denote $T_{\alpha_1}\ldots T_{\alpha_n}.$  For the empty word $\emptyset$
we define $|\emptyset| = 0$ and $T_\emptyset = I$.

The
unitary $U: \m H \otimes \m K \to \m H \otimes \m P$ from our model can be decomposed as
\begin{equation}
  U(h \otimes \ki) = \sum_{j = 1}^d E_j^*h \otimes \epsilon_j \mbox{~~for~} h \in \m
  H,
\end{equation}
where $E_j$'s are some operators in $\m B(\m H),$  for $j =
1,\ldots,d.$  Likewise there exist some operators $C_j$'s in $\m B(\tilde {\m H})$ such that
\begin{equation}
\tilde U(\tilde h \otimes \ki) = \sum_{j = 1}^d C_j^*\tilde h
\otimes \epsilon_j \mbox{~~for~} \tilde h \in \tilde {\m
H}.
\end{equation}
Observe
that $\sum_{j=1}^{d}E_jE_j^* = I$ and $\sum_{j=1}^{d}C_jC_j^* =
I,$ i.e., $\uu E$ and $\uu C$ are coisometric tuples. By equation
(1.2)
$$E_j^*\tilde h = C_j^* \tilde h \mbox{~~for all } \tilde h \in \tilde{\m H}, j = 1,\ldots,d.$$
We recall from \cite{DG11} that such tuple $\uu E = (E_1,\ldots
,E_d)$ is called a \textit{lifting} of $\uu C = (C_1,\ldots
,C_d).$

From a physicist perspective our model is a Markovian approximation of the repeated interaction between a quantum system and a stream of copies of another
quantum system in such a way that there is no backaction. The change of an observable $X \in \m B(\m H)$ until time $n,$ compressed to $\m H,$ is written as
\[Z_n(X):= P_{\m H} U(n)^* (X \otimes I) U(n)|_{\m H}.\]
From equation (1.3) it follows that $Z_n(X)=Z^n(X)$ where $Z(X)= \\ \sum^d_{i=1} E_i X E^*_i:B(\m H) \to B(\m H)$ and $Z$ is called the {\em transition operator} of
the noncommutative Markov chain.

In section 2 we develop our generalized repeated interaction model and
obtain a coisometric operator which intertwines between the minimal
isometric dilations of $\uu E$ and $\uu C,$ and which will
be  crucial for the further investigation in this article.
Using this an outgoing Cuntz scattering system in the sense of \cite{BV05} is constructed for
our model in section 3. Popescu introduced the minimal isometric dilation in \cite{Po89a} and the characteristic
function in \cite{Po89b} of a row contraction, and systematically developed an extensive theory of row contractions (cf. \cite{Po99}, \cite{Po06}).
We use some of the concepts from Popescu's theory in this work.

For the outgoing Cuntz scattering system
in section 4 we give a $\tilde \Lambda$-linear system with an input-output formalism.
A multi-analytic operator appears here as the transfer function and
in the next section we show that this transfer function can be derived from the intertwining coisometry of section 2. In the
scattering interpretation of the transfer function this now mediates between two processes.
This together with a nice product formula obtained in Proposition 2.1 tells us that this identification of transfer function is a reminiscent of the scattering
operator construction using wave operators in Lax-Phillips scattering theory \cite{LP67}, equation (1.5) (cf. \cite{RS79}), with one of the
processes moving forward combined with the other moving backward.
In \cite{YK03} and \cite{GGY08} there are other approaches to
transfer functions. Several works on transfer functions and on
quantum systems using linear system theory can be found in recent
theoretical physics and control theory surveys.

In section 5 we investigate in regard to our model what the notion of observability
implies for the scattering theory and the theory of liftings.
Some techniques used here are similar to those of scattering theory of
noncommutative Markov chains introduced in \cite{KM00}. Characteristic
functions for liftings, introduced in \cite{DG11}, are multi-analytic operators
which classify certain class of liftings. Our model generalizes the setting of \cite{Go10}, and a comparison
is done in section 6 between the transfer function of our model and the characteristic function for the associated lifting
using the series expansion of the transfer function obtained in section 4. As a consequence mathematically generalized interaction
models get firmly linked into the theory of functional models.

\section{A Generalised Repeated Interaction Model}

We begin with three Hilbert spaces $\m
H, \m K$ and $\m P$ with unit vectors $\ki \in \m K$ and $\oi \in \m P$, and unitaries $U$ and $\tilde U$ as in
equation (1.2). In $\m K_\infty = \bigotimes_{i=1}^\infty  \m K$ and $ \m
P_\infty = \bigotimes_{i=1}^\infty  \m P$  define
$\ki_\infty := \bigotimes_{i=1}^\infty \ki$ and $\oi_\infty :=
\bigotimes_{i=1}^\infty  \oi$ respectively. We denote $m$-th
copy of $\ki$ in $\ki_\infty$ by $\ki_m$
 and in terms of this we introduce the notation
$\Omega_{[m,n]}^{{\mathcal K}} := \Omega_m^{\mathcal K}
\otimes \cdots \otimes \Omega_n^{\mathcal K}$. Identify $\m K_{[m,n]}$
with $\Omega_{[1,m-1]}^{\mathcal K}\otimes \m K_{[m,n]}\otimes
\Omega_{[n+1,\infty)}^{\mathcal K},$ $\m H$ with $\m H \otimes
\ki_\infty$ as a subspace of $\m H \otimes \m K_\infty$ and
$\tilde{\m H}$ with $\tilde{\m H} \otimes \ki_\infty$ as a
subspace of $\tilde{\m H} \otimes \m K_\infty$. Similar notations with respect to $\m P$
are also used.
 For simplicity we assume that $d$ is finite but all the results here can be derived also for $d = \infty$.

Associate a row contraction $\underline E$ to the unitary $U$ as in equation (1.3) and define isometries
\[
\widehat V_j^E(h\otimes \eta)  := U^*(h \otimes \epsilon_j)\otimes \eta ~\mbox{for}~j = 1,\ldots ,d,
\]
on the elementary tensors $h\otimes \eta \in \m H \otimes \m
K_\infty$ and extend it linearly to obtain $ \widehat V_j^E \in \m
B(\m H \otimes \m K_\infty)~ \mbox{for}~j =1,\ldots,d$.
We recall that a lifting $\uu T= (T_1, \ldots, T_d)$ of any row contraction $\uu S=(S_1, \ldots, S_d)$ is called its isometric dilation if $T_i$'s are isometries with orthogonal ranges. It can be easily verified that $ \underline {\widehat  V}^E = (\widehat V_1^E,\ldots
,\widehat V_d^E)$ on the space $\m H \otimes \m K_\infty$ is an
 isometric dilation of $\underline E = (E_1,\ldots,E_d).$
If $h \in \m H$ and $k_1 \in \m K,$ then there exist $ h_i \in \m H$ for $i=1,\ldots,d$ such that $ U^*(\sum^d_{i=1} h_i \otimes \epsilon_i)=h \otimes k_1$
 because $U$ is a unitary. This implies
\[ \sum^d_{i=1}  \widehat V_i^E(h_i \otimes  \ki_\infty) = h \otimes k_1 \otimes  \ki_{[2,\infty)} .\]
In addition if $k_2 \in \m K,$ then
\[ \sum^d_{i=1}  \widehat V_i^E(h_i \otimes k_2 \otimes  \ki_{[2,\infty)}) =  U^*(\sum^d_{i=1} h_i \otimes \epsilon_i) \otimes k_2  \otimes \ki_{[3,\infty)} = h \otimes k_1 \otimes k_2  \otimes  \ki_{[3,\infty)}.\]
By induction we conclude that
\[ \m H \otimes \m K_\infty= \overline{span} \{ \widehat V_\alpha^E(h \otimes  \ki_\infty): h \in \m H, \alpha \in \tilde \Lambda \}, \]
i.e., $\underline {\widehat  V}^E$ is the minimal isometric dilation of
$\underline E.$ Note that the minimal isometric dilation is unique up to unitary equivalence (cf. \cite{Po89a}).

Similarly, associate a row contraction $\underline C$ to the unitary $\tilde{U}$ as in equation (1.4) and define isometries
\begin{equation}
 \widehat V_j^C(\tilde h\otimes \eta) := \tilde U^*(\tilde h \otimes \epsilon_j)\otimes \eta ~\mbox{for}~j = 1,\ldots ,d
\end{equation}
on the elementary tensors $ \tilde h \otimes \eta \in \tilde{\m H} \otimes \m K_\infty$ and extend it linearly to obtain
$\widehat V_j^C \in \m B(\tilde{\m H} \otimes \m K_\infty)~ \mbox{for}~j =1,\ldots ,d.$ The tuple $\underline {\widehat V}^C = (\widehat V_1^C,\ldots ,\widehat V_d^C)$ on
 the space $\tilde{\m H} \otimes \m K_\infty$ is the minimal  isometric dilation of $\underline C = (C_1,\ldots,C_d)$. Recall that
\[
 U_m : \mathcal H \otimes \m K_\infty \to \mathcal H \otimes \m K_{[1,m-1]} \otimes \m P_m \otimes \m K_{[m+1, \infty)}
\]
is nothing but  the operator which acts  as $U$ on $\m H \otimes \m
K_m$ and fixes other factors of the infinite tensor products.
Similarly, we define $\tilde U_m$ using
$\tilde{U}.$
\begin{prop} Let $P_n := P_{\tilde {\m H}}
\otimes  I_{\m P_{[1, n]}} \otimes  I_{\m K_{[n+1, \infty )}} \in \m B(\m H \otimes \m P_{[1,n]} \otimes \m K_{[n+1, \infty )})$ for
$ n \in  \mathbb N$. Then
\[
 sot-\displaystyle\lim_{n\to\infty} \tilde U_1^*\ldots \tilde U_n^*P_n U_n\ldots U_1
\]
exists and this limit defines a coisometry $\widehat W :\m  H \otimes \m  K_\infty \to \tilde{\m H} \otimes \m  K_\infty$. Its adjoint $\widehat W^* : \tilde{\m H} \otimes \m  K_\infty \to \m  H \otimes \m  K_\infty$ is given by
\[
  \widehat W^* = sot-\displaystyle\lim_{n\to\infty} U_1^*\ldots U_n^*\tilde U_n\ldots\tilde U_1.
\]
Here sot stands for the strong operator topology.
\end{prop}
\begin{proof} At first we construct the adjoint $\widehat W^*$.  For that consider the dense subset
$\bigcup_{m \geq 1} \tilde{\m H} \otimes \m K_{[1, m]}$ of
$\tilde{\m H} \otimes \m  K_\infty$ and let an arbitrary
 simple tensor element of this dense subset be $\tilde h\otimes k_1 \otimes \ldots \otimes k_\ell \otimes \ki_{[\ell+1, \infty)}$ for
some $\ell \in \mathbb N, \tilde h \in \tilde {\m H}$ and $k_i \in \m K_i.$
Set $a_p = U_1^*\ldots U_p^*\tilde U_p\ldots \tilde U_1(\tilde h \otimes k_1 \otimes \ldots \otimes k_\ell \otimes \ki_{[\ell+1, \infty)})
 ~\mbox{for}~ p \in \mathbb N$. Since $U(\tilde h \otimes \ki) = \tilde U(\tilde h \otimes \ki)~\mbox{for all }~\tilde h
 \in \tilde{\m H}$, we have $a_\ell = a_{\ell+n}~\mbox{for all}~ n\in \mathbb N$. Therefore we deduce that
 \[
 \lim_{n\to\infty} U_1^*\ldots U_n^*\tilde U_n\ldots \tilde U_1(\tilde h \otimes k_1 \otimes \ldots \otimes k_\ell \otimes \ki_{[\ell+1, \infty)})
 \]
exists. Because $U$ and $\tilde U$ are unitaries, we obtain an
isometric extension $\widehat W^*$ to the whole of $\tilde {\m H}
\otimes \m K_\infty$.
 Thus its adjoint is a coisometry $\widehat W :\m  H \otimes \m  K_\infty
\to \tilde{\m H} \otimes \m  K_\infty$.\\
\indent Now we will derive the limit form for $\widehat W$ as claimed in the statement of the proposition.
If $h\otimes \eta \in \m H \otimes \m K_{[1,k]}, \tilde h
\otimes \tilde {\eta}
 \in \tilde{\m H} \otimes \m K_{[1,n]}$ and $k \leq n,$ then
\begin{eqnarray*}
 \langle \widehat W (h\otimes \eta), \tilde h \otimes \tilde{\eta}\rangle & = & \langle h \otimes \tilde {\eta}, \widehat W^*(\tilde h \otimes \tilde {\eta})\rangle\\
 & = & \langle h \otimes  {\eta},  U_1^*\ldots U_n^*\tilde U_n \ldots \tilde U_1(\tilde h \otimes \tilde {\eta})\rangle\\
  & = & \langle \tilde U_1^*\ldots \tilde U_n^*P_n U_n \ldots U_1(h \otimes  {\eta}), \tilde h  \otimes \tilde {\eta}\rangle.
\end{eqnarray*}
Consequently  $\widehat W =  sot-\displaystyle\lim_{n\to\infty} \tilde U_1^*\ldots \tilde U_n^*P_n U_n\ldots U_1$ on a dense subset and therefore
 it can be extended to the whole of
$\m H \otimes \m K_\infty.$
\end{proof}

 Observe that
\begin{equation}
\widehat W^*(\tilde h \otimes \ki_\infty) = \tilde h
\otimes \ki_\infty \mbox{~for all~} \tilde h \in \tilde {\m H}.
\end{equation}
Next we show that this coisometry $\widehat W$ intertwines between $\widehat V_j^E$ and $\widehat V_j^C$ for all $j = 1,\ldots ,d$. For $j= 1,\ldots ,d$, define
\begin{eqnarray*}
 S_j : \m H \otimes \m K_\infty & \to & \m H \otimes \m P_1 \otimes \m K_{[2, \infty )},\\
 h \otimes \eta & \mapsto & h\otimes \epsilon_j \otimes \eta.
\end{eqnarray*}
The following are immediate:
\begin{itemize}
 \item [(1)] $S_j^*(h \otimes p_1 \otimes \eta)  = \langle \epsilon_j, p_1\rangle (h \otimes \eta)~\mbox{for}~(h \otimes p_1 \otimes \eta) \in \m H \otimes \m P_1 \otimes  \m K_{[2, \infty)}$.
 \item [(2)]  $\widehat V_j^E(h\otimes \eta) = U_1^*S_j(h \otimes \eta)~\mbox{for}~ h\otimes \eta \in \m H \otimes \m K_\infty$.
 \item [(3)]  $\widehat V_j^C(\tilde h \otimes \eta) = \tilde U_1^*S_j(\tilde h \otimes \eta)~\mbox{for}~ \tilde h \otimes \eta \in \tilde {\m H}
\otimes \m K_\infty$.
\end{itemize}
\begin{prop}
If $\widehat W$ is as in Proposition 2.1, then
 $$\widehat W \widehat V_j^E = \widehat V_j^C\widehat W, \hspace{.5cm} \widehat V_j^E \widehat W^* = \widehat W^* \widehat V_j^C \mbox{~for all~} j = 1,\ldots, d.$$
\end{prop}
\begin{proof}
If $h \in \m H,  \eta \in \m K_\infty, \tilde h \in \tilde {\m H}$ and  $ k_i \in  \m K_i,$ then by the three observations that were noted preceding this proposition
we obtain for $j = 1,\ldots,d$
 \begin{eqnarray*}
 & &\langle~ \widehat W \widehat V_j^E (h \otimes \eta),~ \tilde h\otimes k_1 \otimes \ldots \otimes k_\ell \otimes \ki_{[\ell+1, \infty)}~\rangle\\
 &= &\langle ~U^*(h \otimes \epsilon_j)\otimes  \eta ,~U_1^*\ldots U_\ell^*\tilde U_\ell \ldots \tilde U_1(\tilde h \otimes k_1 \otimes \ldots
\otimes k_\ell \otimes \ki_{[\ell+1, \infty)})~\rangle.
 \end{eqnarray*}
Substituting $\tilde U(\tilde h \otimes k_1)= \sum_i \tilde h^{(i)}\otimes k_1^{(i)} $ where $\tilde h^{(i)}  \in \tilde {\m H}$ and $ k_1^{(i)} \in \m K$ we obtain
 \begin{eqnarray*}
& &\langle~ \widehat W \widehat V_j^E (h \otimes \eta),~ \tilde h\otimes k_1 \otimes \ldots \otimes k_\ell \otimes \ki_{[\ell+1, \infty)}~\rangle\\
&=&\langle ~h \otimes \epsilon_j \otimes  \eta,~ U_2^* \ldots U_\ell^*\tilde U_\ell \ldots \tilde U_2(\sum_i (\tilde h^{(i)}\otimes k_1^{(i)})\\
&& \otimes k_2\otimes\ldots \otimes k_\ell \otimes \ki_{[\ell+1, \infty)})~\rangle\\
 &=& \sum_i\langle \epsilon_j,~k_1^{(i)}\rangle~\langle ~h \otimes \eta, ~\widehat W^*(\tilde h^{(i)} \otimes k_2
\otimes \ldots \otimes k_\ell \otimes \ki_{[\ell+1, \infty)})~\rangle\\
 &=& \langle~ \widehat W(h \otimes \eta), ~S_j^*\tilde U_1(\tilde h\otimes k_1 \otimes \ldots\otimes k_\ell \otimes \ki_{[\ell+1, \infty)})~\rangle\\
 &=&\langle~ \tilde U_1^*S_j \widehat W((h \otimes \eta),~\tilde h\otimes k_1 \otimes \ldots \otimes k_\ell \otimes \ki_{[\ell+1, \infty)}~\rangle\\
 &=& \langle ~\widehat V_j^C \widehat W(h \otimes \eta),~\tilde h\otimes k_1 \otimes \ldots \otimes k_\ell \otimes \ki_{[\ell+1, \infty)}~\rangle.
 \end{eqnarray*}
Hence $\widehat W\widehat V_j^E = \widehat V_j^C\widehat W$ for all
$j = 1,\ldots ,d.$ To obtain the other equation of the proposition we again use the last
two of the three observations as follows: For $j = 1,\ldots,d$
\begin{eqnarray*}
&&\widehat W^* \widehat V_j^C( \tilde h \otimes k_1 \otimes \ldots \otimes k_\ell \otimes  \ki_{[\ell+1, \infty)})\\
&=& \widehat W^* \tilde U_1^* ( \tilde h\otimes \epsilon_j \otimes k_1 \otimes \ldots \otimes k_\ell \otimes  \ki_{[\ell+2, \infty)})\\
&=& U_1^*U_2^*\ldots U_{\ell+1}^*\tilde U_{\ell+1} \ldots \tilde U_2 \tilde U_1 \tilde U_1^*
( \tilde h\otimes \epsilon_j \otimes k_1 \otimes \ldots \otimes k_\ell \otimes  \ki_{[\ell+2, \infty)})\\
&=& U_1^*U_2^*\ldots U_{\ell+1}^*\tilde U_{\ell+1} \ldots \tilde U_2 S_j ( \tilde h \otimes k_1 \otimes \ldots \otimes k_\ell \otimes  \ki_{[\ell+1, \infty)})\\
&=& U_1^* S_j U_1^* \ldots U_\ell^* \tilde U_\ell \ldots \tilde U_1( \tilde h \otimes k_1 \otimes \ldots \otimes k_\ell \otimes  \ki_{[\ell+1, \infty)})\\
&=& \widehat V_j^E \widehat W^*( \tilde h \otimes k_1 \otimes \ldots \otimes k_\ell \otimes  \ki_{[\ell+1, \infty)})
\end{eqnarray*}
\end{proof}

 Further define
\[
(\m H \otimes  \m K_\infty )^{\circ} := (\m H \otimes  \m K_\infty
) \ominus (\tilde {\m H} \otimes \ki_\infty),
\]
\begin{equation}
 (\tilde {\m H} \otimes  \m K_\infty )^{\circ} := (\tilde {\m H} \otimes  \m K_\infty ) \ominus (\tilde {\m H} \otimes \ki_\infty)\mbox{
and}~ \m H^\circ := \m H \ominus \tilde {\m H}.
\end{equation}
Let $\sum^k_{i=1} \xi_i \otimes \eta_i \in (\m H \otimes \m K_\infty)^\circ$ and $\tilde h \in
\tilde {\m H}.$ Then for $j = 1,\ldots, d$
\begin{eqnarray*}
\langle  \widehat V_j^E (\sum_i \xi_i \otimes \eta_i ), \tilde h \otimes \Omega^{\m K}_\infty \rangle
&=& \langle\sum_i U^* (\xi_i \otimes \epsilon_j) \otimes \eta_i, \tilde h \otimes \Omega^{\m K}_\infty \rangle\\
&=& \langle  \sum_i \xi_i \otimes \epsilon_j \otimes \eta_i, \tilde U (\tilde h \otimes \Omega^{\m K}_1) \otimes \Omega^{\m K}_{[2, \infty)} \rangle =0
\end{eqnarray*}
because $\tilde U$ maps into $\tilde {\m H} \otimes \m P$ and $\sum^k_{i=1} \xi_i \otimes \eta_i \bot \tilde{\m H} \otimes \Omega^{\m K}.$
Therefore $\widehat V_j^E (\m H \otimes  \m K_\infty
)^{\circ} \subset (\m H \otimes  \m K_\infty )^{\circ}$ for $j =
1, \ldots, d.$ Similarly
$\widehat V_j^C (\tilde{\m H} \otimes  \m K_\infty )^{\circ}
\subset \m (\tilde{\m H} \otimes  \m K_\infty )^{\circ}$ for $j =
1,\ldots, d.$ Set $V_j^E := \widehat V_j^E|_{(\m H \otimes \m
K_\infty )^{\circ}}$ and $V_j^C := \widehat V_j^C|_{(\tilde{\m H}
\otimes \m K_\infty )^{\circ}}$ for $j = 1,\ldots, d.$ If we define
$$W^* := \widehat W^*|_{(\tilde {\m H}
\otimes \m K_\infty )^{\circ}},$$
then by equation (2.2) it follows that $W^* \in \m B((\tilde {\m H}
\otimes  \m K_\infty )^{\circ}, (\m H \otimes  \m K_\infty
)^{\circ})$. The operator $W^*$ is an isometry because it is a
restriction of an isometry and $W$, the adjoint of $W^*$, is the
restriction of $\widehat W$ to $(\m H \otimes  \m K_\infty
)^{\circ},$ i.e., $W = \widehat W|_{( \m H \otimes  \m K_\infty
)^{\circ}}$.
\begin{rem}
It follows that
\[
 W V_j^E =  V_j^C W
\]
for $j =1,\ldots,d$.
\end{rem}

\section{Outgoing Cuntz Scattering Systems}

In this section we aim to construct an outgoing Cuntz scattering
system (cf. \cite {BV05}) for our model. This will assist us in the next section to work with an input-output formalism and
to associate a transfer function to the model.

Following are some notions from the multivariable operator theory:

\begin{defn} Suppose $\underline T = (T_1,\ldots, T_d)$ is a row contraction where $T_i \in \m B (\m L).$
\begin{itemize}
\item [(1)] If $T_i$'s are isometries with orthogonal ranges, then
the tuple $\underline T = (T_1,\ldots,\\ T_d)$ is called a
\textit{row isometry}.

\item [(2)] If $ \overline{span}_{j=1,\ldots,d} T_j \m L = \m L $
and $\underline T = (T_1,\ldots,T_d)$ is a row isometry, then
$\underline T$ is called a  \textit{row unitary}.

\item [(3)] If there exist a subspace $\m E$ of $\m L$ such that
$\m L = \bigoplus_{\alpha \in \tilde {\Lambda}} T_\alpha \m E$ and
$\underline T = (T_1,\ldots,T_d)$ is a row isometry, then
  $\uu T$ is called a \textit{row shift} and $\m E$ is called a \textit{wandering subspace} of $\m L$ w.r.t. $\uu T.$
\end{itemize}
\end{defn}
\begin{defn} A collection $(\m L, \underline V = (V_1,\ldots, V_d), \m G_{*}^+, \m G )$ is called an outgoing Cuntz scattering system
(cf. \cite{BV05}), if  $\underline V$ is a row isometry on the Hilbert space $\m L$, and $\m G_{*}^+$ and $ \m G$ are subspaces of
$\m L$ such that
\begin{itemize}
 \item [(1)] for $\m E_* := \m L \ominus \overline{span}_{j =1,\ldots , d} V_j \m L$, the tuple $\underline V|\m G_{*}^+$ is a row shift where $\m G_{*}^+ = \bigoplus _{\alpha \in \tilde \Lambda}  V_\alpha \m E_*$.
\item [(2)] there exist  $\m E := \m G \ominus \overline{span}_{j =1,\ldots, d} V_j \m G$ with $\m G =  \bigoplus _{\alpha \in \tilde \Lambda}  V_\alpha \m E$, i.e., $\underline V|_\m G$ is a row shift.
\end{itemize}
\end{defn}

In the above definition the part (1) is the Wold decomposition (cf. \cite{Po89a}) of the row isometry $\underline V$ and therefore $\m G_{*}^+$ can
be derived from $\underline V.$ But $\m G_{*}^+$ is included in the data because it helps in describing the scattering phenomenon.
We continue using the notations from the previous section.  $\widehat V_j^E$'s are isometries with orthogonal ranges and
because  $(\epsilon_j)_{j = 1}^{d}$ is an orthonormal basis of $\m
P,$ we have
\[
\overline {span}_{j =1,\ldots,d} \widehat V_j^E(\m H \otimes \m K_\infty) = \m H \otimes \m K_\infty.
\]
Thus $\widehat {\uu V}^E$ is a row unitary on $\m H \otimes  \m K_\infty $. Now using the fact that $V^E_j =
\widehat V^E_j|_{(\m H \otimes \m K_\infty)^\circ} $ we infer that $V^E_j$'s are isometries with orthogonal ranges.
Therefore $\uu V^E$ is a row isometry on $(\m H \otimes \m K_\infty)^\circ$ .

\begin{prop}
If  $\m Y :=  \tilde {\m H}  \otimes (\ki_1)^\perp \otimes \ki_{[2,\infty)} \subset \tilde{\m H} \otimes \m K_\infty$, then
\[
 W^*\m Y \perp \overline {span}_{j =1,\ldots,d} \ V_j^E(\m H \otimes \m K_\infty)^{\circ}.
\]
\end{prop}
\begin{proof}
 By Proposition 2.1 it is easy to see that
\begin{equation}
W^* \m Y = U^*_1 \tilde {U}_1 \m Y \subset \m H \otimes \m K_1 \otimes \ki_{[2, \infty)}.
\end{equation}
 Let $\tilde h_i \in \tilde {\m H}$ and $k_i \perp \ki_1$ for $i =1,\ldots,n,$ i.e.,  $\textstyle \sum_i \tilde h_i \otimes k_i \otimes \ki_{[2, \infty)} \in \m Y.$ For  $\sum_k h_k \otimes \eta_k \in (\m H \otimes \m K_\infty)^{\circ}$ with $h_k \in \m H$ and $\eta_k \in \m K_\infty$
\begin{eqnarray*}
&&\langle~ W^* (\textstyle \sum_i \tilde {h_i} \otimes k_i \otimes \ki_{[2, \infty)}), V_j^{E}(\sum_k h_k \otimes \eta_k)~\rangle\\
 & = & \langle ~U^*\tilde {U}(\textstyle \sum_i \tilde {h_i} \otimes k_i) \otimes \ki_{[2, \infty)}, \sum_k U^*(h_k \otimes \epsilon_j )\otimes \eta_k~\rangle \\
& = & \langle ~\tilde {U}(\textstyle \sum_i \tilde {h_i} \otimes k_i) \otimes \ki_{[2, \infty)}, \sum_k h_k \otimes \epsilon_j \otimes \eta_k~\rangle
= 0.
\end{eqnarray*}
The last equality holds because $\sum_k h_k \otimes \eta_k  \perp \tilde {\m H} \otimes \ki _\infty$. Thus $W^*\m Y \perp \overline {span}_{j =1,\ldots,d}  V_j^E(\m H \otimes \m K_\infty)^{\circ}.$
\end{proof}

The following Proposition gives an explicit description of the Wold decomposition of $\uu V^E:$
\begin{prop}
If $\m Y$ is defined as in the previous proposition, then $W^*\m Y$ is a wandering subspace  of $\uu V^E$, i.e., $V_\alpha^E(W^* \m Y) \perp V_\beta^E(W^* \m Y)$ whenever $\alpha, \beta \in \tilde{\Lambda}$, $\alpha \neq \beta$, and
\[
 W^*\m Y = (\m H \otimes  \m K_\infty )^{\circ} \ominus \overline {span}_{j =1,\ldots,d} V_j^E(\m H \otimes  \m K_\infty )^{\circ}.
 \]
\end{prop}
\begin{proof}
 By Proposition  3.3 it is immediate that $V_\alpha^E(W^* \m Y) \perp V_\beta^E(W^* \m Y)$ whenever $\alpha, \beta \in \tilde{\Lambda}$, $\alpha \neq \beta$ and
$W^*\m Y \subset (\m H \otimes  \m K_\infty )^{\circ} \ominus
\overline {span}_{j =1,\ldots,d} V_j^E(\m H \otimes  \m K_\infty )^{\circ}$. The only thing that remains to be shown is that
\[
(\m H \otimes  \m K_\infty )^{\circ} \ominus
\overline {span}_{j =1,\ldots,d} V_j^E(\m H \otimes  \m K_\infty )^{\circ} \subset W^*\m Y.
\]
\indent Let  $x \in (\m H \otimes  \m K_\infty )^{\circ} \ominus
\overline{span}_{j =1,\ldots ,d} V_j^E(\m H \otimes  \m K_\infty
)^{\circ}$. Write down the decomposition of $x$ as $x_1 \oplus
x_2$ w.r.t.  $W^*\m Y \oplus (W^*\m Y)^{\perp}$. So
 $x - x_1 = x_2 $ is orthogonal to both $\overline{span}_{j =1,\ldots ,d} V_j^E(\m H \otimes  \m K_\infty )^{\circ}$ and $W^*\m Y$. Now we show that
if any element in $(\m H \otimes  \m K_\infty )^{\circ}$ is orthogonal to $\overline {span}_{j =1,\ldots ,d} V_j^E(\m H \otimes  \m K_\infty )^{\circ}$
and $W^*\m Y$, then it is the zero vector. Let $x_0$ be such an element. Because $x_0 \in (\m H \otimes  \m K_\infty )^{\circ}$ and $x_0 \perp W^*\m Y,$
\[
x_0 \perp U^*(\tilde {\m H} \otimes \epsilon_j)\otimes \ki_{[2,\infty)}
\]
for $j=1,\ldots,d$. This implies $x_0 \perp \overline {span}_{j
=1,\ldots ,d}\widehat V_j^E (\tilde{\m H}\otimes \ki_\infty).$
We also know that
 \[x_0 \perp \overline {span}_{j =1,\ldots,d} V_j^E (\m H \otimes  \m K_\infty )^{\circ} (= \overline {span}_{j =1,\ldots,d}\widehat V_j^E (\m H \otimes  \m K_\infty )^{\circ}).
 \]
Therefore
\[
 x_0 \perp \overline {span}_{j =1,\ldots,d}\widehat V_j^E (\m H \otimes  \m K_\infty ).
\]
Since $\widehat {\uu V}^E$ is a row unitrary, $x_0 \perp \m H \otimes \m K_\infty$. So $x_0 = 0$ and hence $x = x_1 \in W^* \m Y.$ We conclude that $(\m H \otimes  \m K_\infty )^{\circ} \ominus
\overline {span}_{j =1,\ldots,d} V_j^E(\m H \otimes  \m K_\infty )^{\circ} \subset W^*\m Y$.
\end{proof}
\begin{prop}
 If $\m E := \m H \otimes (\ki_1)^\perp \otimes \ki_{[2, \infty)} \subset (\m H \otimes  \m K_\infty )^{\circ}$, then $V^E_\alpha \m E \perp V^E_\beta \m E$ whenever $\alpha, \beta \in \tilde{\Lambda}, \alpha \neq \beta$ and $ (\m H \otimes  \m K_\infty )^\circ =  \m H^\circ \oplus \m \bigoplus_{\alpha \in \tilde \Lambda}V^E_{\alpha} \m E$.
\end{prop}

\begin{proof}
If $|\alpha| = |\beta|$ and $\alpha \neq \beta$, then it is easy to see that $V^E_\alpha \m E \perp V^E_\beta \m E$ because
ranges of $V_j^E$'s are mutually orthogonal. If $|\alpha| \neq |\beta|$ (without loss of generality  we can assume that $|\alpha| > |\beta|$), then by taking the inner product at the tensor factor $\m K_{|\alpha |+1}$ we obtain $V^E_\alpha \m E \perp V^E_\beta \m E$.\\
\indent To prove the second part of the proposition, observe that for $n\in \mathbb{N}$,
\begin{eqnarray*}
 && \m H \otimes \m K_{[1, n]} \otimes \ki_{[n+1, \infty]} \\ 
 &=& (\m H \otimes \ki_\infty ) \oplus (\m H \otimes
 (\ki_1)^\perp \otimes \ki_{[2, \infty)} )\oplus (\m H \otimes \m K_1\otimes \\
 && (\ki_2)^\perp \otimes \ki_{[3, \infty)}) \oplus \cdots
 \oplus ~(\m H \otimes \m K_{[1, n-1]} \otimes (\ki_n)^\perp \otimes \ki_{[n+1, \infty)})\\
 &=&  (\tilde {\m H} \otimes \ki_\infty )\oplus (\m H^\circ \otimes \ki_\infty) \oplus \m E \oplus \displaystyle\bigoplus_{j=1}^{d} V_j^E \m E \oplus \cdots \oplus \displaystyle\bigoplus_{|\alpha |= n-1}^{d} V_\alpha ^E \m E.
\end{eqnarray*}
Taking $n \to \infty$ we have the following:
\[
 \m H \otimes \m K_\infty = (\tilde {\m H} \otimes \ki_\infty )\oplus (\m H^\circ \otimes \ki_\infty) \oplus \displaystyle\bigoplus_{\alpha \in \tilde \Lambda} V_\alpha ^E \m E.
 \]
Since $(\m H \otimes \m K_\infty)^{\circ} = (\m H \otimes \m K_\infty) \ominus (\tilde {\m H} \otimes \ki_\infty )$, it follows that
 \[
(\m H \otimes  \m K_\infty )^\circ = \m H^\circ \oplus \displaystyle\bigoplus_{\alpha \in \tilde \Lambda} V_\alpha ^E \m E.
 \]
\end{proof}
\indent We sum up Propositions 3.3, 3.4 and 3.5 in the
following theorem:
\begin{thm}
For a generalized repeated interaction model involving unitaries $U$ and $\tilde U$ as before set
$\m Y :=  \tilde {\m H}  \otimes (\ki_1)^\perp \otimes \ki_{[2,\infty)}$ and
$\m E := \m H \otimes (\ki_1)^\perp \otimes \ki_{[2, \infty)}.$
 If $\m E_* := W^* \m Y$, $\m G^+_* := \bigoplus_{\alpha \in \tilde \Lambda} V_\alpha ^E \m E_*$ and $ \m G := \bigoplus_{\alpha \in \tilde \Lambda} V_\alpha ^E \m E$,
then the collection
$$((\m H \otimes  \m K_\infty )^{\circ}, {\uu V}^E = (V_1^E,\ldots, V_d^E), \m G_{*}^+, \m G )$$
is an outgoing Cuntz scattering system such that  $(\m H \otimes  \m K_\infty )^\circ = \m H^\circ \oplus  \m G $.
\end{thm}

\begin{rem}
\indent Applying arguments similar to those used for proving the
second part of the  Proposition 3.5 one can prove the following:
  \[
   (\tilde {\m H} \otimes \m K_\infty)^\circ = \bigoplus_{\alpha \in \tilde \Lambda}V^C_{\alpha} \m Y.
  \]
\end{rem}
We refer the reader to Proposition 3.1 of \cite{Go11} for a result in a similar direction.

\section{$\tilde \Lambda$-Linear Systems and Transfer Functions}

We would demonstrate that the outgoing Cuntz scattering system $((\m H
\otimes \m K_\infty )^{\circ}, {\uu V}^E = (V_1^E,\ldots, V_d^E),
\m G_{*}^+, \m G )$ from Theorem 3.6 has interesting relations with
a generalization of the linear systems theory that is associated to our interaction model.
For a given model involving unitaries $U$ and $\tilde U$ as before, let us define the input
space as
\[
\m U  :=  \m E = \m H \otimes (\ki_1)^\perp \otimes \ki_{[2, \infty)} \subset (\m H \otimes  \m K_\infty )^{\circ}
\]
and the output space as
\[
\m Y  =  \tilde {\m H} \otimes (\ki_1)^\perp \otimes \ki_{[2, \infty)} \subset (\tilde {\m H} \otimes \m K_\infty)^\circ.
\]

Here we assume that a quantum system $\m A$ interacts with a stream of copies of another quantum system $\m B$ and we assume $\m H$ is the (quantum mechanical)
Hilbert space of $\m A.$ Let $\m K_i$ be the Hilbert space of a part of a stream of copies of
 $\m B$  at time $i$ immediately before the interaction with $\m A.$ Let the Hilbert space $\m P_i$ be that the part of a stream of copies of $\m B$ at time $i$
immediately after the interaction with $\m A.$  $\ki$ and $\oi$ denote states indicating that no copy of quantum system $\m B$
is present and so no interaction is taking place at time $i.$ Then $\eta \in \m U= \m H \otimes (\ki_1)^\bot \otimes \ki_{[2,\infty)} \subset \m H \otimes
\m K_\infty$ represents a vector state with copies of quantum system $\m B$ arriving at time 1 and stimulating an interaction between the stream of copies of $\m A$ and
$\m B,$ but no further copy of $\m B$ arriving at later times. But some activity is induced which goes on for a longer period.

 Note that $\m H \otimes \m K = \m H \oplus \m U$ and $\tilde{\m H} \otimes \m K = \tilde{\m H} \oplus \m Y$.
So $U$ maps $\m H \oplus \m U$ onto  $\m H \otimes \m P$ and $\tilde U$ maps $\tilde{\m H}\oplus \m Y$ onto  $\tilde{\m H} \otimes \m P$.
Using unitaries $U$ and $\tilde U$ we  define
$ F_j : \m H \to \m U ~~\mbox {and}~~ D_j:\tilde {\m H} \to \m Y~~ \mbox{for}~ j =1,\ldots,d$ by
\begin{align}
 \displaystyle \sum_{j=1}^d F_j^* \eta \otimes \epsilon_j  :=  U (0 \oplus \eta),~~~
\displaystyle \sum_{j=1}^d D_j^* y\otimes \epsilon_j :=  \tilde U(0 \oplus y) \mbox{~~for~}\eta \in \m U~\mbox{and}~y\in \m Y.
\end{align}
  Combining equation (4.1) with equations (1.3) and (1.4) we have for $h \in \m H,~ \eta \in \m U,~\tilde h \in \tilde{\m H}~\mbox{and}~y\in \m Y$
\begin{eqnarray}
 U (h \oplus \eta) = \displaystyle \sum_{j=1}^d(E_j^* h + F_j^* \eta)\otimes
 \epsilon_j,\\
 \tilde U(\tilde h \oplus y)  = \displaystyle \sum_{j=1}^d (C_j^* \tilde h +D_j^* y)\otimes \epsilon_j
\end{eqnarray}
respectively.
Using equation (4.3) it can be checked that
\begin{eqnarray}
  \tilde U^* (\tilde h \otimes \epsilon_j) =  C_j \tilde h \oplus D_j\tilde h
\mbox{~~for~} \tilde h \in \tilde {\m H}; j=1,\ldots,d.
\end{eqnarray}

 Let us define
$$\tilde C := \sum_{j=1}^d D_jP_{\tilde{\m H}} E_j^*: \m H \to \m Y,~~~\tilde D := \sum_{j=1}^d D_jP_{\tilde{\m H}} F_j^*: \m U \to \m Y$$
where $P_{\tilde{\m H}}$ is the orthogonal projection onto $\tilde{\m H}$. It follows that
\begin{eqnarray}
 P_{\m Y}\tilde U^*P_1U(h \oplus \eta) = \tilde C h +\tilde D \eta
\end{eqnarray}
where $h \in \m H, \eta \in \m U, P_1$ is as in Proposition 2.1 and $P_{\m Y}$ is the orthogonal
projection onto $\m Y$.\\
  Define a \textit{colligation} of operators (cf. \cite{BV05}) using the operators $E_j^*$'s, $F_j^*$'s, $\tilde C~\mbox{and}~ \tilde D$ by
\[
 \m C_{U,\tilde U} :=  \begin{pmatrix}
          E_1^* & F_1^* \\
      \vdots& \vdots\\
        E_d^* & F_d^*\\
     \tilde C & \tilde D
         \end{pmatrix}
: \m H \oplus \m U \to \displaystyle\bigoplus_{j = 1}^d \m H \oplus \m Y.
\]
From the colligation $\m C_{U,\tilde U}$ we get the following $\tilde \Lambda$-\textit{linear system}
$\sum_{U,\tilde U}$:
\begin{eqnarray}
 x(j\alpha) &=& E_j^* x(\alpha) + F_j^* u(\alpha), \\
~ y(\alpha) &=&  \tilde Cx(\alpha) + \tilde D u(\alpha)
\end{eqnarray}
where $j = 1,\ldots,d$ and $\alpha ,j\alpha$ are words in $ \tilde \Lambda$, and
\[
  x : \tilde \Lambda \to \m H,~~~ u : \tilde \Lambda \to \m U,~~~ y :\tilde \Lambda \to \m Y.
\]
If $x(\emptyset)$ and $u$ are known, then using $\sum_{U,\tilde
U}$ we can compute $x$ and $y$ recursively. Such a $\tilde \Lambda$-linear
system is also called a noncommutative Fornasini-Marchesini system
in \cite{BGM06} in reference to \cite{FM78}.

Let $z = (z_1 ,\ldots, z_d )$ be a $d$-tuple of formal
noncommuting indeterminates. Define the Fourier transforms of $x,u$
and $y$ as
\[
 \hat x(z) = \displaystyle \sum_{\alpha \in \tilde \Lambda} x(\alpha) z^{\alpha}, ~~~
\hat u(z) =  \sum_{\alpha \in \tilde \Lambda} u(\alpha) z^{\alpha}, ~~~\hat y(z) = \sum_{\alpha \in \tilde \Lambda} y(\alpha) z^{\alpha}
\]
respectively where $z^\alpha=z_{\alpha_n}\ldots  z_{\alpha_1}$ for
$\alpha = \alpha_n\ldots \alpha_1 \in \tilde \Lambda$.  Assuming that $z$-variables commute with the coefficients the input-output relation
\[
 \hat y(z) = \Theta_{U,\tilde U} (z) \hat u(z)
\]
can be obtained on setting $x(\emptyset): = 0$  where
\begin{equation}
 \Theta_{U,\tilde U}(z) := \displaystyle \sum_{\alpha \in \tilde \Lambda} \Theta_{U,\tilde U}^{(\alpha)}
z^\alpha := \tilde D + \tilde C  \displaystyle \sum_{\beta \in \tilde \Lambda, j = 1,\ldots,d}(E_{\bar\beta})^* F_j^* z^{\beta j}.
 \end{equation}
Here $\bar \beta=\beta_1 \ldots \beta_n$ is the reverse of $\beta=\beta_n \ldots \beta_1
\in \tilde \Lambda$ and $\Theta_{U,\tilde U}^{(\alpha)}$ maps $\m U$ to $\m Y$. The
formal noncommutative power series $\Theta_{U,\tilde U}$ is called
the \textit{transfer function} associated to the unitaries $U$ and
$\tilde U$. The transfer function is a  mathematical tool for
encoding the evolution of a $\tilde \Lambda$-linear system.
For $y(\alpha) \in \m Y$ with
$\sum_{\alpha \in \tilde \Lambda}\|y(\alpha)\|^2 < \infty$,
any series $\sum_{\alpha \in \tilde \Lambda} y(\alpha) z^\alpha$ stands for a series converging to an element of $\ell^2(\tilde \Lambda, \m Y)$.

\begin{thm}
 The map $ M_{\Theta_{U,\tilde U}} : \ell^2 (\tilde \Lambda, \m U) \to  \ell^2 (\tilde \Lambda, \m Y)$ defined by
 \[
 M_{\Theta_{U,\tilde U}} \hat u(z) := \Theta_{U,\tilde U}(z)\hat u(z)
 \]
 is a contraction.
\end{thm}
\begin{proof}
Observe that $P_{\m Y}\tilde U^*P_1U(\tilde h \otimes \ki_\infty)
= 0$ for all $\tilde h \in  \tilde {\m H}$. Consider another
colligation which is defined as follows:
\[
 \m C^\circ_{U,\tilde U} := \begin{pmatrix}
          E_1^{*\circ} & F_1^{*\circ} \\
      \vdots& \vdots\\
        E_d^{*\circ}& F_d^{*\circ}\\
     \tilde C^\circ & \tilde D
         \end{pmatrix}
: \m H^\circ \oplus \m U \to \displaystyle\bigoplus_{j = 1}^d \m H^\circ \oplus \m Y
\]
where $E_j^{*\circ}:= P_{\m H^\circ}E_j^*|_{\m H^\circ} : \m
H^\circ \to  \m H^\circ,~F_j^{*\circ} := P_{\m H^\circ} F_j^* : \m
U \to  \m H^\circ$ and $\tilde C^\circ := \tilde C|_{\m H^\circ}
:\m H^\circ \to  \m Y$ for $j=1,\ldots,d$. Recall that $\m H^\circ$ and $(\m H \otimes  \m K_\infty )^{\circ}$ were defined in equation array (2.3). Consider the outgoing
Cuntz scattering system $((\m H \otimes  \m K_\infty )^{\circ},
{\uu V}^E = (V_1^E,\ldots, V_d^E), \m G_{*}^+, \m G )$, with $(\m
H \otimes  \m K_\infty )^\circ = \m H^\circ \oplus  \m G $,
constructed by us  in Theorem 3.6. In Chapter 5.2 of \cite{BV05} it is shown that there is an associated  unitary colligation
\begin{eqnarray}
 \begin{pmatrix}
          \hat E_1& \hat F_1 \\
      \vdots& \vdots\\
        \hat E_d & \hat F_d\\
    \hat M & \hat N
         \end{pmatrix}
: \m H^\circ \oplus \m E \to  \displaystyle\bigoplus_{j = 1}^d \m H^\circ \oplus \m E_*
\end{eqnarray}
such that $(\hat E_j, \hat F_j) = P_{\m H^\circ} (V_j^E)^*|_{\m H^\circ \oplus \m E}$ and $(\hat M, \hat N) = P_{\m E_*}|_{\m H^\circ \oplus \m E}.$
Recall that $\m E$ and $\m E_*$ were introduced in Proposition 3.5 and Theorem 3.6 respectively. From equations (4.2) and (4.5) we observe that
$(E_j^{*\circ}, F_j^{*\circ}) = P_{\m H^\circ \otimes
\epsilon_j}U|_{\m H^\circ \oplus \m E}$ (identifying $\m H^\circ $
with $\m H^\circ \otimes \epsilon_j$) and $(\tilde C^\circ, \tilde
D) = P_{\m Y}\tilde U^*P_1U|_{\m H^\circ \oplus \m E}$. Using
these observations we obtain the following relations:
\begin{eqnarray}
 U^*(E_j^{*\circ}, F_j^{*\circ}) &=& U^*P_{\m H^\circ \otimes \epsilon_j} U|_{\m H^\circ \oplus \m E} = P_{U^*(\m H^\circ \otimes \epsilon_j)}
|_{\m H^\circ \oplus \m E} = P_{V_j^E \m H^\circ}|_{\m H^\circ \oplus \m E} \notag \\
 &=& V_j^E P_{\m H^\circ} (V_j^E)^*|_{\m H^\circ \oplus \m E} = V^E_j (\hat E_j, \hat F_j)
\end{eqnarray}
for $j=1,\ldots,d$ and
\begin{eqnarray}
 U^*\tilde U(\tilde C^\circ, \tilde D) &=& U^*\tilde U P_{\m Y} \tilde U^* P_1 U|_{\m H^\circ \oplus \m E}
 =  U^*P_{\tilde U \m Y}P_1U|_{\m H^\circ \oplus \m E} = U^*P_{\tilde U \m Y}U|_{\m H^\circ \oplus \m E} \notag \\
 &=& P_{U^* \tilde U \m Y}|_{\m H^\circ \oplus \m E}
 = P_{W^* \m Y}|_{\m H^\circ \oplus \m E} \hspace{.5cm} (\mbox{by equation (3.1)}) \notag \\
&=& P_{\m E_*}|_{\m H^\circ \oplus \m E} = (\hat M, \hat N).
\end{eqnarray}
\indent Let  $\hat u(z) = \sum_{\alpha \in \tilde \Lambda}
u(\alpha) z^{\alpha}  \in  \ell^2 (\tilde \Lambda, \m U)$ with
$u(\alpha) \in \m U$ such that \\ $\sum_{\alpha \in \tilde \Lambda}\|
u(\alpha)\|^2 < \infty$. We would prove that
\[
 \|M_{\Theta_{U,\tilde U}} \hat u(z)\|^2 \leq \|\hat u(z)\|^2.
\]
 Define $x : \tilde \Lambda \to \m H$ by equation (4.6) such that $x(\emptyset) = 0.$
 Further, define $x^{\circ}(\alpha) := P_{\m H^\circ}x(\alpha)$ for all $\alpha \in \tilde \Lambda.$
 Now applying the projection $P_{\m H^\circ}$ to relation (4.6)
 on both sides
 and using the fact $\tilde {\m H}$ is invariant under $E_j^*$  for $j=1,\ldots,d$ we obtain the following
 relation:
\begin{equation}
 x^\circ(j \alpha) = E_j^{*\circ} x^\circ(\alpha) + F_j^{*\circ} u(\alpha)~\mbox {for all}~\alpha \in \tilde \Lambda,~j=1,\ldots,d.
\end{equation}
Because $P_{\m Y}\tilde U^*P_1U_1(\tilde h \otimes
\ki_\infty) = 0$ for all $\tilde h \in  \tilde {\m H}$ we conclude
by equation (4.5) that
\begin{equation}
 \tilde C\tilde h=0  \mbox{~~for~} \tilde h \in \tilde {\m H}.
\end{equation}
This implies
\begin{equation}
\tilde C x(\alpha) = \tilde C^\circ x^\circ(\alpha) \mbox{~for
all~} \alpha \in \tilde \Lambda.
\end{equation}

Define $y :\tilde \Lambda \to \m Y $ by
\begin{equation}
y(\alpha) :=  \tilde Cx(\alpha) + \tilde D u(\alpha)
\end{equation}
for all $\alpha\in  \tilde \Lambda$.  Recall that the input-output relation stated just before the theorem is
 \[
 \hat y(z)=  \sum_{\alpha \in \tilde \Lambda} y(\alpha) z^{\alpha} = \Theta_{U,\tilde U} (z) \hat u(z)(= M_{\Theta_{U,\tilde U}} \hat u(z)).
 \]
Using the unitary colligation given in equation (4.9) we have
\begin{eqnarray*}
 \|x^\circ(\alpha)\|^2 + \|u(\alpha)\|^2 &=&  \displaystyle \sum_{j=1}^d \|\hat E_jx^\circ(\alpha) + \hat F_ju(\alpha)\|^2+ \|\hat Mx^\circ(\alpha) + \hat N u(\alpha)\|^2\\
  &=& \displaystyle \sum_{j=1}^d \|E_j^{*\circ}x^\circ(\alpha) + F_j^{*\circ}u(\alpha)\|^2 + \|\tilde C^\circ x^\circ(\alpha) + \tilde D u(\alpha)\|^2\\
  &=& \displaystyle \sum_{j=1}^d \|x^\circ(j\alpha)\|^2 + \|\tilde C x(\alpha) + \tilde D u(\alpha)\|^2\\
  &=& \displaystyle \sum_{j=1}^d \|x^\circ(j\alpha)\|^2 +\| y(\alpha)\|^2
\end{eqnarray*}
for all $\alpha \in \tilde \Lambda$. In the above calculation
equations (4.10), (4.11), (4.12), (4.14) and (4.15) respectively have been used. This
gives us
\begin{eqnarray*}
    \|u(\alpha)\|^2 -\| y(\alpha)\|^2 = \displaystyle \sum_{j=1}^d \|x^\circ(j\alpha)\|^2 - \|x^\circ(\alpha)\|^2
\end{eqnarray*}
for all $\alpha \in \tilde \Lambda$. Summing over all $\alpha \in
\tilde \Lambda$ with $|\alpha| \leq n$ and using  the fact that
$x^\circ(\emptyset) = 0$ we obtain
\[
\displaystyle \sum_{|\alpha| \leq n} \|u(\alpha)\|^2 - \displaystyle \sum_{|\alpha| \leq n}\| y(\alpha)\|^2 = \displaystyle \sum_{|\alpha| = n+1}\|x^\circ(\alpha)\|^2 \geq 0~\mbox{for all}~ n\in \mathbb{N}.
\]
Therefore
\[
 \displaystyle \sum_{|\alpha| \leq n}\| y(\alpha)\|^2 \leq \displaystyle \sum_{|\alpha| \leq n} \|u(\alpha)\|^2 ~\mbox{for all}~ n\in \mathbb{N}.
\]
Finally taking limit $n\to\infty$ both the sides we get that $M_{\Theta_{U,\tilde U}}$ is a contraction.
\end{proof}
\indent $ M_{\Theta_{U,\tilde U}}$ is a \textit{multi-analytic
operator} (\cite{Po95}) (also called \textit{analytic
intertwining operator} in \cite{BV05}) because
\[
 M_{\Theta_{U,\tilde U}}( \displaystyle \sum_{\alpha \in \tilde \Lambda} u(\alpha) z^\alpha z^j) = M_{\Theta_{U,\tilde U}}( \displaystyle \sum_{\alpha \in \tilde \Lambda} u(\alpha) z^\alpha )z^j ~~\mbox{for}~ j =1, \ldots, d,
 \]
i.e., $M_{\Theta_{U,\tilde U}}$ intertwines with right
translation. The noncommutative power series $\Theta_{U,\tilde U}$
is called the \textit{symbol} of $M_{\Theta_{U,\tilde U}}$.

\section{Transfer Functions, Observability and Scattering}

We would now establish that the transfer function can be derived from the coisometry $W$ of section 2. In the last section $d$-tuple
$z = (z_1 ,\ldots, z_d )$  of formal noncommuting indeterminates were employed. Treat
 $(z^\alpha)_{\alpha \in \tilde \Lambda}$ as an orthonormal basis of $\ell^2(\tilde \Lambda, \mathbb C).$
Assume $\m Y$ and $\m U$ to be the spaces associated with our model with unitaries $U$ and $\tilde U$ as in the last section.
It follows from Remark 3.7 that there exist a unitary operator $\tilde{\Gamma}:  (\tilde{\m H} \otimes  \m K_\infty )^{\circ} \to  \ell^2 (\tilde \Lambda, \m Y)$ defined by
\[
 \tilde{\Gamma}(V^C_\alpha y) := y z^{\bar \alpha}~ \mbox { for all }~ \alpha \in \tilde \Lambda, y \in \m Y.
\]
We observe the following intertwining relation:
\begin{align}
 \tilde{\Gamma}(V^C_\alpha y) = (\tilde{\Gamma}y)z^{\bar \alpha}.
 \end{align}
Similarly, using Theorem 3.6, we can define a unitary operator $ \Gamma : (\m H \otimes  \m K_\infty )^{\circ} (= (\m H^\circ \oplus \m G)) \to  \m H^\circ \oplus \ell^2(\tilde \Lambda, \m U)$ by
 \[
 \Gamma(\mathring h \oplus V^E_\alpha \eta) := \mathring h \oplus \eta z^{\bar \alpha} \mbox { for all }~ \alpha \in \tilde \Lambda
 \]
where $\mathring h \in \m H^\circ, \eta \in\m  U.$ In this case the intertwining relation is
\begin{align}
 \Gamma (V^E_\alpha \eta) = (\Gamma \eta) z^{\bar \alpha}.
\end{align}
Using the coisometric operator $W$, which appears in Remark 2.3, we define $\Gamma_W$ by the following commutative diagram:
\begin{equation}
\xymatrix{
(\m H \otimes \m K_\infty)^\circ \ar[r]^{{ W}} \ar[d]_{\Gamma}
& (\tilde{\m H} \otimes  \m K_\infty )^{\circ} \ar[d]^{\tilde{\Gamma}}
\\
\m H^{\circ} \oplus \ell^2 (\tilde \Lambda, \m U) \ar[r]^{\Gamma_W}
&   \ell^2 (\tilde \Lambda, \m Y),
}
\end{equation}
i.e., $\Gamma_W  = \tilde{\Gamma} W \Gamma^{-1}$.

\begin{thm}
$\Gamma_W$ defined by the above commutative diagram satisfies
 \[
 \Gamma_W|_{\ell^2 (\tilde \Lambda,~ \m U)} =  M_{\Theta_{U, \tilde U}}.
\]
\end{thm}
\begin{proof}
Using the intertwining relation $V_j^C W = W V_j^E$ from Remark 2.3, and equations (5.1) and (5.2) we obtain
\begin{eqnarray*}
 \Gamma_W(\eta z^\beta z^j) &=& \tilde{\Gamma} W\Gamma^{-1} (\eta z^\beta z^j) = \tilde{\Gamma} WV^E_jV^E_{\bar\beta}\eta\\
 &=& \tilde{\Gamma} V^C_jV^C_{\bar\beta}W\eta = (\tilde{\Gamma} W\eta)z^\beta z^j = \Gamma_W(\eta z^\beta)z^j
\end{eqnarray*}
for $\eta \in \m U, \beta \in \tilde \Lambda, j=1,\ldots,d$. Hence, $ \Gamma_W|_{\ell^2 (\tilde \Lambda,~ \m U)}$ is a multi-analytic operator.
For computing its symbol we determine $\Gamma_W \eta$ for $\eta \in \m U$, where $\eta$ is identified  with $\eta z^\phi \in \ell^2(\tilde \Lambda, \m U).$ For $\alpha = \alpha_{n-1}\ldots \alpha_1 \in \tilde\Lambda$ let $P_\alpha$ be the orthogonal projection onto
\begin{eqnarray*}
 & & \tilde{\Gamma}^{-1}\{f \in \ell^2(\tilde \Lambda, \m Y) : f = yz^\alpha ~\mbox{for some}~ y\in \m Y \} \\
& =& V^C_{\bar \alpha}\m Y = \tilde U_1^*\ldots \tilde U_{n-1}^*(\tilde {\m H} \otimes \epsilon_{\alpha_1} \otimes \cdots \otimes \epsilon_{\alpha_{n-1}}\otimes (\ki_n)^\perp \otimes \ki_{[n+1,\infty)})
\end{eqnarray*}
with $\tilde U_i$'s as  in Proposition 2.1.\\
\indent Recall that the tuple $\underline E$ associated with the unitary $U$ is a lifting of the tuple $\underline C$
(associated with the unitary $\tilde U$) and so $\underline E$
can be written as a block matrix in terms of $\underline C$
 as follows:  $E_j =
\begin{pmatrix}
          C_j & 0 \\
        B_j & A_j\\
         \end{pmatrix}$
for $j=1,\ldots,d$ w.r.t. to the decomposition $\m H = \tilde {\m
H} \oplus \m H^\circ$ where $\uu B$ and $\uu A$ are some row
contractions. Because $\uu E$ is a coisometric lifting of $\uu C$
we have $$\displaystyle \sum_{j=1}^d C_jC_j^* = I \mbox{~and~}
\displaystyle \sum_{j=1}^d C_jB_j^* = 0$$ (cf. \cite{DG11}) . Now
using these relations and equations (4.2), (4.3) and (4.4) it can
be easily verified that
\[
 P_\alpha \tilde U_1^*\ldots \tilde U_{n}^*P_nU_n \ldots U_1\eta = P_\alpha \tilde U_1^* \ldots \tilde U_{m}^*P_mU_m\ldots U_1\eta ~\mbox{for all}~ m \geq n, \eta \in \m U.
\]
Using the formula of $W$ from Proposition 2.1 we obtain
\[
 P_\alpha W\eta  = P_\alpha \tilde U_1^*\ldots \tilde U_{n}^*P_nU_n\ldots U_1\eta~\mbox{for}~  \eta \in \m U.
\]
Finally for $\eta \in \m U$
\begin{eqnarray*}
P_\alpha \tilde U_1^*\ldots \tilde U_{n}^*P_nU_n\ldots U_1\eta = \left\{\begin{array}{ll} \tilde D \eta ~~~~~~~~~~~~~~~~~~~~~~~~~~~~~~~~~~~~~~~~~~~~~~ \mbox{if}~~ n=1,\alpha = \emptyset,\\
V^C_{\bar\alpha}(\tilde C E^*_{\alpha_{n-1}}\ldots E^*_{\alpha_2}F^*_{\alpha_1}\eta) ~~~~~~~~~~~~~~~~~\mbox{if}~~   n=|\alpha |+1 \geq 2.
 \end{array}\right.
\end{eqnarray*}
This implies for $\eta \in \m U$
\[\tilde{\Gamma} W \Gamma^{-1} \eta = \tilde{\Gamma} W \eta =  \tilde D \eta \oplus \sum_{|\alpha| \geq1} (\tilde C E^*_{\alpha_{n-1}}\ldots
E^*_{\alpha_2}F^*_{\alpha_1}\eta) z^\alpha. \]
Comparing this with equation (4.8) we conclude that $ \Gamma_W|_{\ell^2 (\tilde \Lambda,~ \m U)} =  M_{\Theta_{U, \tilde U}}$.
\end{proof}

Note that the Theorem 4.1 and its proof concern the transfer function of the $\tilde{\Lambda}$-linear system and has nothing to do with the scattering
 theory. Theorem 5.1, on the other hand, is the scattering theory part in the sense of Lax-Phillips \cite{LP67}. The same function $M_{\Theta_{U, \tilde U}}$
relates the outgoing Fourier representation
for a vector in the ambient scattering Hilbert space to the incoming Fourier representation for the same vector. This makes $M_{\Theta_{U, \tilde U}}$
the scattering function for the outgoing Cuntz scattering system.
We introduce a notion from the linear systems theory for our model:
\begin{defn}
 The observability operator $W_0: \m H ^\circ \to \ell^2(\tilde \Lambda, \m Y)$ is defined as the restriction of the operator $\Gamma_W$ to $\m H ^\circ $, i.e., $W_0 = \Gamma_W|_{\m H ^\circ }$.
 \end{defn}
\noindent  It follows  that $W_0\mathring h = (\tilde C (E_{\bar\alpha})^*\mathring h)_{\alpha \in \tilde \Lambda}$.
  Popescu has studied the similar types of operators called Poisson kernels in \cite {Po99}.
 \begin{defn}
If there exist $k,K > 0$ such that for all $\mathring h \in \m H ^\circ$
\[
 k\|\mathring h\|^2 \leq \displaystyle \sum_{\alpha \in \tilde \Lambda} \|\tilde C (E_{\bar\alpha})^*\mathring h\|^2 = \|W_0\mathring h\|^2 \leq K\| \mathring h\|^2,
\]
then the  $\tilde \Lambda-$linear system is called (uniformly) observable.
\end{defn}
\noindent We illustrate below that the notion of observability
is closely related to the scattering theory notions of noncommutative Markov chains.
Observability of a system for $\dim \m H < \infty$ is interpreted as the property of the system that in the absence of $\m U$-inputs we can determine the original state $h\in \m H^\circ$ of the system from all $\m Y$-outputs at all times. Uniform observability is an analog of this for $\dim \m H =\infty.$

We extend $W_0$ to
 \[
  \widehat W_0 :(\tilde {\m H} \oplus \m H^\circ)(= \m H) \longrightarrow \tilde {\m H} \oplus \ell^2(\tilde \Lambda, \m Y)
 \]
by defining $\widehat W_0 \tilde h := \tilde h$ for all $\tilde h \in \tilde {\m H}$. If $W_0$ is uniformly observable, then using $\hat k = k$ and $\hat K = \mbox{max}\{1, K\}$ the above inequalities can be extended to $\widehat W_0$ on $\m H$ as
\[
   \hat k\|h\|^2 \leq \  \|\widehat W_0h\|^2 \leq \hat K\| h\|^2
 \]
for all $h \in \m H$.

Before stating the main theorem of this section regarding
observability  we recall from \cite{DG11} the following: Let $\uu
C$ be a row contraction on a Hilbert space $\m H_C.$ The lifting
$\underline E$ of $\underline C$ is called \textit{subisometric}
\cite{DG11} if the minimal isometric dilations $\underline {\widehat V}^E$  and
$\underline {\widehat V}^C$  of $\uu E$ and $\uu C$ respectively are
unitarily equivalent and the corresponding unitary, which
intertwines between $\widehat V_i^E$ and  $\widehat V_i^C$ for all
$i=1, 2, \ldots, d,$ acts as identity on $\m H_C.$ Some of the techniques used here are from the scattering theory of noncommutative
Markov chains (cf. \cite{KM00}, \cite{Go04}).

\begin{thm}
 For any $\tilde \Lambda$-linear system associated to a generalized repeated interaction model with unitaries $U,\tilde U$ the following statements are equivalent:
\begin{itemize}
 \item [(a)] The system is (uniformly) observable.
 \item [(b)] The observability operator $W_0$ is isometric.
\item [(c)] The tuple $\underline E$ associated with the unitary $U$ is a subisometric lifting of the tuple
$\underline C$ (associated with the unitary $\tilde U$).
\item[(d)] $W : (\m H \otimes \m K_\infty)^\circ \to (\tilde {\m H} \otimes \m K_\infty)^\circ$ is unitary.
\end{itemize}
If one of the above holds, then
\begin{itemize}
 \item [(e)] The transfer function $\Theta_{U,\tilde U}$ is inner, i.e., $M_{\Theta_{U,\tilde U}} : \ell^2(\tilde \Lambda, \m U) \to \ell^2(\tilde \Lambda, \m Y)$ is isometric.
\end{itemize}
If we have additional assumptions, viz.  $\dim \m H < \infty$ and
 $\dim \m P \geq 2$, then the converse  holds, i.e., $(e)$ implies all
of $(a), (b), (c)$ and  $(d)$.
\end{thm}
\begin{proof}
 Clearly $(d) \Rightarrow (b) \Rightarrow (a)$. We now prove $(a) \Rightarrow (d)$. Because the system is
 (uniformly) observable there exist $k > 0$
such that for all $\mathring h \in \m H ^\circ$
\[
 k\|\mathring h\|^2 \leq  \|W_0\mathring h\|^2.
\]
Since $\bigcup_{m \geq 1} \m H \otimes \m K_{[1, m]}$ is a dense subspace of  $\m H \otimes \m K_\infty$, for any $0 \neq \eta \in \m H \otimes \m K_\infty$
there exist $n \in \mathbb N$ and $\eta^\prime \in \m H \otimes \m K_{[1, n]}$ such that
 \[
  \| \eta - \eta^\prime \| < \dfrac{\sqrt{k}}{\sqrt{k}+1}\|\eta\|.
 \]
Let $\eta_0 \in \m H \otimes \m K_{[1, n]}$. Suppose $U_n \ldots U_1 \eta_0 = h_0 \otimes p_0 \otimes \ki_{[n+1, \infty)}$, where $h_0 \in \m H$, $p_0 \in \m P_{[1,n]}.$ Then clearly
 $$\displaystyle\lim_{N\to\infty}\|\tilde U^*_1\ldots \tilde U^*_n\tilde U^*_{n+1} \ldots \tilde U^*_NP_NU_N \ldots U_{n+1}U_n\ldots U_1 \eta_0 \|  = \|\widehat W_0 h_0\| \|p_0\|$$
and thus by Proposition 2.1 it is equal to $\|\widehat W
\eta_0\|$. Because the system is (uniformly) observable,
$$\|\widehat W_0 h_0\| \|p_0\| \geq \sqrt{k}\|h_0\|\|p_0\|.$$
Therefore $\|\widehat W \eta_0\|^2 \geq k
\|\eta_0\|^2.$ However, in general $U_n\ldots U_1 \eta_0 = \sum_j
h_0^{(j)} \otimes p_0^{(j)} \otimes \ki_{[n+1, \infty)}$ with
$h_0^{(j)} \in \m H$ and some mutually orthogonal vectors $p_0^{(j)}
\in \m P_{[1,n]}$. By using the above inequality for each term of
the summation and then adding them we find that
 in general for all $\eta_0 \in \m H \otimes
\m K_{[1, n]}$
\[
 \|\widehat W \eta_0\|^2 \geq k\|\eta_0\|^2.
\]
In particular, for $\eta^\prime \in \m H \otimes \m K_{[1, n]}$ we have the above inequality. Therefore
\begin{eqnarray*}
 \|\widehat W \eta \| &\geq & \|\widehat W \eta^\prime\| - \|\widehat W(\eta^\prime - \eta) \|\\
 & \geq & \sqrt{k} \|\eta^\prime \|-\|\eta- \eta^\prime\|\\
  & \geq & \sqrt{k}\|\eta\|- (\sqrt{k}+1)\|\eta- \eta^\prime\| > 0.
 \end{eqnarray*}
This implies $\widehat W \eta \neq 0$ for all $0 \neq \eta \in \m H
\otimes \m K_\infty$ and hence $\widehat W$ is injective. Recall
that $\widehat W$ is a coisometry  and an injective coisometry is
unitary. Further, because $\widehat W(\tilde h \otimes \ki_\infty)
= \tilde h \otimes \ki_\infty$ for all
$\tilde h \in \tilde {\m H}$ it follows that $W$ is unitary. This establishes $(a)  \Rightarrow (d)$ and we have proved $(a) \Leftrightarrow (b) \Leftrightarrow (d)$.\\
\indent Next we prove $(d) \Leftrightarrow (c)$. Assume that $(d)$ holds. Since  $W$ is unitary, clearly $\widehat W$ is unitary. We know that $\widehat W$ intertwines between the minimal isometric dilations $\widehat {\uu V}^E$ and $\widehat
{\uu V}^C$ of $\uu E$ and $\uu C$ respectively. Hence $\uu E$ is a subisometric lifting of $\uu C$.\\
\indent Conversely, if we assume $(c)$, then by the definition of subisometric lifting there exist a unitary operator $$\widehat W_1 : \m H \otimes \m K_\infty \longrightarrow
\tilde{\m H }\otimes \m K_\infty$$ which intertwines between $\widehat {\uu V}^E$ and $\widehat {\uu V}^C,$ and $\widehat W_1 $ acts as an identity on
$\tilde {\m H} \otimes \ki_\infty.$ To prove $W$ is unitary it is enough to prove $\widehat W$ is unitary. We show that $\widehat W = \widehat W_1$.
By the definition of the minimal isometric dilation we know that $\tilde{\m H }\otimes \m K_\infty = \overline {span}\{ \widehat V_\alpha^C (\tilde h \otimes \ki_\infty)
: \tilde h \in \tilde {\m H},\alpha \in \tilde \Lambda\}$.
For $j = 1,\ldots,d$ and $\tilde h \in \tilde {\m H},$ by equation (2.2) and Proposition 2.2,
\begin{eqnarray*}
 \widehat W^*\widehat V_j^C (\tilde h \otimes \ki_\infty) &=& \widehat V_j^E\widehat W^* (\tilde h \otimes \ki_\infty) = \widehat V_j^E (\tilde h \otimes \ki_\infty)\\
 &=& \widehat{W}_1^*\widehat V_j^C \widehat{W}_1(\tilde h \otimes \ki_\infty) = \widehat{W}_1^*\widehat V_j^C(\tilde h \otimes \ki_\infty).
\end{eqnarray*}
Thus $\widehat W^* = \widehat W_1^*$ and hence $\widehat W = \widehat W_1$. \\
\indent To prove $(d) \Rightarrow (e)$ we at first note that since $W$ is unitary, $\Gamma_W$ is also unitary. By Theorem 4.2,  we have
$M_{\Theta_{U,\tilde U}} = \Gamma_W|_{\ell^2(\tilde \Lambda, \m U)}$. Since a restriction of a unitary operator is an isometry, $M_{\Theta_{U,\tilde U}}$ is isometric. \\
\indent Finally with the additional assumptions $\dim \m H < \infty$
and $\dim\m P \geq 2$, we show $(e) \Rightarrow (b)$.  Define
\[
 \m H_{scat} := \m H \cap\widehat W^*(\tilde {\m H} \otimes \m K_\infty) = \tilde{\m H} \oplus \{\mathring{h} \in \m H^\circ : \|W_0 \mathring h\| = \|\mathring h\|\}.
\]
 Since $\|\widehat W_0 h\| = \displaystyle\lim_{n\to\infty}\|\tilde U_1\ldots \tilde U_n\tilde P_nU_n\ldots U_1 h \|$ by  Proposition 2.1, the following can be easily verified:
 \begin{equation}
  U(\m H_{scat} \otimes \ki) \subset \m H_{scat} \otimes \m P.
 \end{equation}
Because $M_{\Theta_{U,\tilde U}} = \Gamma_W|_{\ell^2(\tilde \Lambda, \m U)}$ is isometric by (e), it can be checked that
\begin{equation}
 U( \m H \otimes (\ki)^\perp) \subset \m H_{scat} \otimes \m P.
\end{equation}
 Combining equations (5.4) and (5.5) we have
 \begin{eqnarray*}
  U^*((\m H \ominus \m H_{scat})\otimes\m P) \subset (\m H \ominus \m H_{scat}) \otimes \ki.
 \end{eqnarray*}
 Since $\dim \m H < \infty$ and $\dim \m P \geq 2$,  we obtain $\m H \ominus \m H_{scat} =
\{0\},$ i.e., $\m H = \m H_{scat}$. This implies $W_0$ is isometric and hence $(e) \Rightarrow (b)$.
 \end{proof}

 \section{Transfer Functions and Characteristic Functions of Liftings}

Continuing with the study of our generalized repeated interaction model,
 from equations (2.1) and (4.4)
  we obtain
\begin{eqnarray}
\widehat V_j^C (\tilde h \otimes  \ki_\infty) = (C_j \tilde h \oplus D_j\tilde h)\otimes \ki_{[2, \infty)} ~~\mbox{for}~ \tilde h \in \tilde {\m H}~\mbox{and}~j=1,\ldots,d.
\end{eqnarray}
Let  $D_C:=(I-\uu C^* \uu C)^\frac{1}{2}:
\bigoplus^d_{i=1} \tilde {\m H} \rightarrow \bigoplus^d_{i=1} \tilde {\m H}$ denote the defect operator and $\m D_C:= \overline{\mbox{Range~} D_C}$. The full Fock space over
$\mathbb{C}^d$ ($d\geq 2$) denoted by $\m F$ is
$$\m F =\mathbb{C}\oplus \mathbb{C}^d \oplus (\mathbb{C}^d)^{\otimes ^2}\oplus \cdots
\oplus (\mathbb{C}^d) ^{\otimes ^m}\oplus \cdots.
$$
The vector $e_\emptyset:=1\oplus 0\oplus \cdots$ is called the vacuum vector.
Let $\{e_1, \ldots , e_d \}$ be the standard
orthonormal basis of $\mathbb{C}^d$.
For $\alpha \in \tilde \Lambda$ and $|\alpha|=n,$ $e_{\alpha}$ denote the vector $e_{\alpha
_1}\otimes e_{\alpha _2}\otimes \cdots \otimes e_{\alpha _n}$ in
the full Fock space $\m F$. We recall that Popescu's construction \cite{Po89a} of the minimal isometric
dilation $\uu{\tilde V}^C = (\tilde V_1^C,\ldots , \tilde V_d^C)$
on $\tilde {\m H} \oplus (\m F \otimes \m D_C)$ of the tuple $\uu C$ is
\[
  \tilde V_j^C(\tilde h \oplus \displaystyle \sum_{\alpha \in \tilde{\Lambda}} e_\alpha \otimes d_\alpha) = C_j\tilde h \oplus [e_\emptyset \otimes (D_C)_j \tilde h + e_j
\otimes  \displaystyle \sum_{\alpha \in \tilde{\Lambda}} e_\alpha \otimes d_\alpha]
\]
for $\tilde h \in \tilde {\m H}$ and $d_\alpha \in \m D_C$
where $(D_C)_j\tilde h = D_C(0,\ldots,\tilde  h, \ldots, 0)$
 ($\tilde h$ is embedded at the $j^{th}$ component).  So
\begin{eqnarray}
 \tilde V_j^C\tilde h = C_j\tilde h \oplus (e_\emptyset \otimes (D_C)_j \tilde h) ~~\mbox{for}~ \tilde h \in \tilde {\m H}~\mbox{and}~ j= 1,\ldots, d.
\end{eqnarray}
From equations (6.1) and (6.2) it follows that
\begin{eqnarray}
 \|\sum_{j = 1}^d D_j\tilde h_j\|^2 = \|\sum_{j = 1}^d (D_C)_j\tilde h_j\|^2
\end{eqnarray}
where $\tilde h_j \in \tilde {\m H}$ for $j= 1,\ldots, d$. Let
$\Phi_C :\overline{span}\{D_j\tilde h : \tilde h \in \tilde{\m H}, j = 1, \ldots, d\} \to \m D_C$ be the unitary given by
\[
 \Phi_C(\sum_{j = 1}^d  D_j\tilde h_j) = \sum_{j = 1}^d  (D_C)_j\tilde h_j ~~~\mbox{for}~ \tilde h_j \in \tilde {\m H}~\mbox{and}~ j= 1,\ldots, d.
\]
Similarly for $E_i$'s and $F_i$'s obtained from interaction $U$ in equation (4.2) we set  $D_E:=(I-\uu E^* \uu E)^\frac{1}{2}:
\bigoplus^d_{i=1} \m H \rightarrow \bigoplus^d_{i=1} \m H$ and $\m D_E := \overline{\mbox{Range~} D_E},$  and define another unitary operator
$  \Phi_E :\overline{span}\{F_j h : h \in \m H, j = 1, \ldots, d\} \to \m D_E$
by
\[
 \Phi_E(\sum_{j = 1}^d F_j h_j) = \sum_{j = 1}^d (D_E)_j h_j ~~\mbox{for}~ h_j \in \m H~\mbox{and}~ j= 1,\ldots, d.
\]

The second equation of (4.1) yields
 \[
  \displaystyle \sum_{j=1}^d D_j D_j^*y = y ~~\mbox{for~}~y \in \m Y.
 \]
This implies
 \[
\overline{span}\{D_j\tilde h : \tilde h \in \tilde{\m H}, j = 1, \ldots, d\} = \m Y.
\]
Similarly, we can show that
$ \overline{span}\{F_j h : h \in \m H, j = 1, \ldots, d\} = \m U.$
Thus $\Phi_C$ is a unitary from $\m Y$ onto $\m D_C$ and $\Phi_E$
is a unitary from $\m U$ onto $\m D_E$.
As a consequence we have for $i,j = 1,\ldots,d$
\begin{eqnarray}
 D_j^*D_i = (D_C)_j^*(D_C)_i = \delta_{ij}I-C_j^*C_i,\\
 F_j^*F_i = (D_E)_j^*(D_E)_i = \delta_{ij}I-E_j^*E_i.
\end{eqnarray}
Define unitaries $\tilde{M}_{\Phi_C} : \ell^2(\tilde \Lambda, \m Y) \to \m F \otimes \m D_C$ and $ \tilde{\Phi}_E : \m U z^\emptyset \to e_\emptyset \otimes \m D_E$ by
\begin{eqnarray*}
 \tilde{M}_{\Phi_C}\big(\displaystyle \sum_{\alpha \in \tilde{\Lambda}} y_\alpha z^\alpha\big)&:=& \displaystyle \sum_{\alpha \in \tilde{\Lambda}} e_
{\bar \alpha} \otimes \Phi_C (y_\alpha),\\
\tilde{\Phi}_E ( u z^\emptyset) &:=& e_\emptyset \otimes \Phi_E u
\end{eqnarray*}
which would be useful in comparing transfer functions with characteristic functions.

 Define $D_{*,A}:= (I-\uu A \uu A^*)^\frac{1}{2}:
\m H^\circ \rightarrow \m H^\circ$  and $\m D_{*,A}:=\overline{\mbox{Range~} D_{*,A}}.$ Because $\uu E$ is a coisometric lifting of $\uu C,$
using Theorem 2.1 of  \cite{DG11} we conclude that there
exist an isometry $\gamma : \m D_{*, A} \to \m D_C$ with
 $\gamma D_{*, A} h = \uu B^*h$ for all $h \in {\m H}^\circ$.
Further, for $h \in {\m H}^ \circ$
\begin{eqnarray*}
\Phi_C \tilde C h
 &=& \Phi_C \displaystyle \sum_{j = 1}^d D_jP_{\tilde{\m H}}E_j^* h
= \Phi_C \displaystyle \sum_{j = 1}^d D_jP_{\tilde{\m H}}(B_j^* h\oplus A_j^* h)\\
&=& \Phi_C \displaystyle \sum_{j = 1}^d D_j B_j^* h
= \displaystyle \sum_{j = 1}^d  (D_C)_j B_j^* h\\
&=& D_C \uu B^*h
= \uu B^*h.
\end{eqnarray*}
The last equality holds because for the coisometric tuple $\uu C$ the operator $D_C$ is  the projection onto
$\m D_C$ and $\mbox{Range~} \uu B^* \subset \m D_C.$ This implies
\begin{equation}
  \Phi_C \tilde C h= \gamma D_{*, A} h.
\end{equation}

The characteristic function $M_{C,E}: \m F \otimes {\m D}_E \to \m F \otimes {\m D}_C$ of lifting $\uu E$ of $\uu C,$ which was introduced in \cite{DG11},
and its symbol $\Theta_{C,E}$ has the following expansion:
For $i = 1,\ldots,d$ and $h \in \tilde{\m H}$ 	
\begin{align}
\Theta_{C,E}  (D_E)_i h = e_\emptyset \otimes
[(D_C)_i h - \gamma D_{*,A} B_i h] - \sum_{|\alpha|\geq 1} e_\alpha \otimes \gamma D_{*,A}
(A_\alpha)^* B_i h,
\end{align}
and for $h \in \m H^\circ$
\begin{eqnarray}
\Theta_{C,E}\, (D_E)_i h &=& - e_\emptyset \otimes \gamma D_{*,A} A_i h \nonumber \\ 
&+& \sum^d_{j=1} e_j \otimes \sum_\alpha e_\alpha \otimes \gamma D_{*,A} (A_\alpha)^*
(\delta_{ji} I -  A^*_j A_i) h.
\end{eqnarray}

\begin{thm} Let $U$ and $\tilde U$ be unitaries associated with a generalized repeated interaction model, and the lifting $\uu E$ of $\uu C$ be the corresponding lifting.
Then the characteristic function $M_{C,E}$ coincides with the transfer function $\Theta_{U, \tilde U},$ i.e.,
\[
 \tilde{M}_{\Phi_C} \Theta_{U, \tilde U}(z)= \Theta_{C,E} \tilde{\Phi}_E.
\]
\end{thm}

\begin{proof}
If $h\in \m H$ and $i=1, \ldots, d,$ then by equation (4.8)
\begin{eqnarray}
 && \tilde{M}_{\Phi_C} \Theta_{U, \tilde U}(z) (F_i  h z^\emptyset)
\nonumber\\ & = &  \tilde{M}_{\Phi_C} [\tilde D ~z^\emptyset+   \displaystyle \sum_{\beta \in \tilde \Lambda, j = 1,\ldots ,d} \tilde C(E_{\bar\beta})^* F_j^* z^{\beta j}] (F_i  hz^\emptyset) \nonumber \\
  & = &  \tilde{M}_{\Phi_C} [\tilde D  F_i  h ~z^\emptyset+ \displaystyle \sum_{\beta \in \tilde \Lambda, j = 1,\ldots ,d} \tilde C  (E_{\bar\beta})^* F_j^*F_i  h~z^{\beta j}].
\end{eqnarray}

 Case 1. $h \in \tilde{\m H}:$
\begin{eqnarray*}
 \tilde D  F_i  h
  &=& \displaystyle \sum_{j = 1}^d D_jP_{\tilde{\m H}}F_j^* F_i h
  = \displaystyle \sum_{j = 1}^d D_jP_{\tilde{\m H}}( \delta_{ij}I-E_j^*E_i) h\\
  &=& D_i h - \big (\displaystyle \sum_{j = 1}^d D_j P_{\tilde{\m H}}E_j^*\big)E_i h
  =  D_i h - \tilde C E_i h\\
  &=&  D_i h - \tilde C (C_i h\oplus B_i h)
   =  D_i h-\tilde C B_i h.
   \end{eqnarray*}
Second and last equalities follows from equations (6.5) and (4.13) respectively. By equation (6.5) again we obtain

\begin{eqnarray*}
&& \displaystyle \sum_{\beta \in \tilde \Lambda, j = 1,\ldots,d}\tilde C (E_{\bar\beta})^* F_j^*F_i h~z^{\beta j}\\
& = &
\displaystyle   \sum_{\beta \in \tilde \Lambda, j = 1,\ldots,d}\tilde C (E_{\bar\beta})^*( \delta_{ij}I-E_j^*E_i) h ~z^{\beta j} \\
&=& \displaystyle \sum_{\beta \in \tilde \Lambda} \tilde C (E_{\bar\beta})^* h~z^{\beta i} - \displaystyle   \sum_{\beta \in \tilde \Lambda, j = 1,\ldots,d}\tilde C (E_{\bar\beta})^* E_j^*E_i h ~z^{\beta j}\\
&=& - \displaystyle   \sum_{\beta \in \tilde \Lambda, j = 1,\ldots,d}\tilde C (E_{\bar\beta})^* E_j^*E_i h ~z^{\beta j}\\
&&(\mbox{because~} \tilde C  (E_{\bar\beta})^*h = \tilde C (C_{\bar\beta})^*h = 0 \mbox{~by equation (4.13)})\\
&=& - \displaystyle   \sum_{\beta \in \tilde \Lambda, j = 1,\ldots,d}\tilde C (E_{\bar\beta})^* \big((C_j^*C_i +B_j^*B_i)h \oplus A_j^*B_ih\big) ~z^{\beta j}\\
&=&  - \displaystyle   \sum_{\beta \in \tilde \Lambda, j = 1,\ldots,d}\tilde C (A_{\bar\beta})^*A_j^*B_ih ~z^{\beta j}\hspace{.5cm}(\mbox{by equation (4.13)})\\
&=&- \displaystyle   \sum_{ |\alpha | \geq 1} \tilde C (A_{\bar\alpha})^* B_ih ~z^{\alpha}.
\end{eqnarray*}
So by equation (6.9) we have for all $i= 1, \ldots, d$ and  $h\in \tilde{\m H}$
\begin{eqnarray*}
 && \tilde{M}_{\Phi_C} \Theta_{U, \tilde U}(z) (F_i h z^\emptyset)\\
 &=&  \tilde{M}_{\Phi_C} [(D_i h-\tilde C B_i h)~z^\emptyset - \displaystyle   \sum_{ |\alpha | \geq 1} \tilde C (A_{\bar\alpha})^* B_ih ~z^{\alpha}]\\
 &=& e_\emptyset \otimes \Phi_C (D_i h-\tilde C B_i h) - \displaystyle   \sum_{ |\alpha | \geq 1} e_{\bar \alpha} \otimes \Phi_C (\tilde C (A_{\bar\alpha})^* B_ih)\\
 &=& e_\emptyset \otimes [(D_C)_ih - \gamma D_{*, A}B_i h] - \displaystyle   \sum_{ |\alpha | \geq 1} e_{\bar \alpha} \otimes \gamma D_{*, A} (A_{\bar\alpha})^* B_ih.
\end{eqnarray*}
By equation (6.7) it follows that
\begin{eqnarray*}
 \tilde{M}_{\Phi_C} \Theta_{U, \tilde U}(z) (F_i h z^\emptyset) &=&  \Theta_{C, E} (e_\emptyset \otimes (D_E)_ih)\\
&=& \Theta_{C, E} \tilde{\Phi}_E (F_i h z^\emptyset).
\end{eqnarray*}

\noindent Case 2. $h \in {\m H}^\circ:$
\begin{eqnarray*}
 \tilde D  F_i  h
 &=& \displaystyle \sum_{j = 1}^d D_jP_{\tilde{\m H}}F_j^* F_i h
 = \displaystyle \sum_{j = 1}^d D_jP_{\tilde{\m H}}( \delta_{ij}I-E_j^*E_i) h\\
 &=& D_i P_{\tilde{\m H}}  h - \big (\displaystyle \sum_{j = 1}^d D_j P_{\tilde{\m H}}E_j^*\big)E_i h
 =  - \tilde C A_i h
  \end{eqnarray*}
Second equality follows from equation (6.5). By equations (6.5) and (4.13) again we obtain
\begin{eqnarray*}
 \displaystyle \sum_{\beta \in \tilde \Lambda, j = 1,\ldots,d}\tilde C (E_{\bar\beta})^* F_j^*F_i h~z^{\beta j}
& = &
\displaystyle   \sum_{\beta \in \tilde \Lambda, j = 1,\ldots,d}\tilde C (E_{\bar\beta})^*( \delta_{ij}I-E_j^*E_i) h ~z^{\beta j} \\
&=&  \displaystyle   \sum_{\beta \in \tilde \Lambda, j = 1,\ldots ,d}\tilde C (A_{\bar\beta})^*(\delta_{ij}I- A_j^*A_i) h~z^{\beta j}.
\end{eqnarray*}
So by equation (6.9) we have for all $i= 1, \ldots, d$ and  $h\in {\m H}^\circ$
\begin{eqnarray*}
&& \tilde{M}_{\Phi_C} \Theta_{U, \tilde U}(z) (F_i h z^\emptyset)\\
&=&  \tilde{M}_{\Phi_C} [ - \tilde C A_i h ~z^\emptyset +  \displaystyle   \sum_{\beta \in \tilde \Lambda, j = 1,\ldots,d}\tilde C (A_{\bar\beta})^*(\delta_{ij}I- A_j^*A_i) h~z^{\beta j}]\\
 &=& - e_\emptyset \otimes \Phi_C(\tilde C A_i h) + \displaystyle   \sum_{\beta \in \tilde \Lambda, j = 1,\ldots ,d} e_j \otimes e_{\bar \beta}  \otimes \Phi_C(\tilde C (A_{\bar\beta})^*(\delta_{ij}I- A_j^*A_i) h)\\
 &=& - e_\emptyset \otimes  \gamma D_{*, A}A_i h + \displaystyle   \sum_{\beta \in \tilde \Lambda, j = 1,\ldots,d} e_j \otimes e_{\bar \beta}  \otimes \gamma D_{*, A} (A_{\bar\beta})^*(\delta_{ij}I- A_j^*A_i) h.
 \end{eqnarray*}
By equation (6.8) it follows that
\begin{eqnarray*}
 \tilde{M}_{\Phi_C} \Theta_{U, \tilde U}(z) (F_i h z^\emptyset) &=&  \Theta_{C, E} (e_\emptyset \otimes (D_E)_ih)\\
&=& \Theta_{C, E} \tilde{\Phi}_E ( F_i h z^\emptyset).
\end{eqnarray*}

Hence we conclude that
\[
 \tilde{M}_{\Phi_C} \Theta_{U, \tilde U}(z)= \Theta_{C,E} \tilde{\Phi}_E.
\]
\end{proof}

The transfer function is a notion affiliated with the input/state/output linear system, while the scattering function is a notion affiliated
with the scattering theory in the sense of Lax-Phillips.
For our repeated interaction model Theorem 6.1 elucidates that the transfer function is identifiable with the characteristic function of the associated lifting. This establishes a
strong connection between a model for quantum systems and the multivariate operator theory. Connections between them were also endorsed in other works
like \cite{Bh96}, \cite{Go04}, \cite{DG07} and \cite{Go11}, and this indicates that such approaches to quantum systems using multi-analytic operators are promising.


\subsection*{Acknowledgment}
The first author received a support from UKIERI to visit Aberystwyth University, UK in July 2011 which was helpful for this project.

\begin{thebibliography}{1}

\bibitem{BGM06}
J. A. Ball, G. Groenewald, T. Malakorn,
\textit{Conservative structured noncommutative multidimensional linear systems. The state space method generalizations and applications},
179--223, Oper. Theory Adv. Appl., 161, Birkh\"{a}user, Basel (2006).


\bibitem{Bh96}
B. V. R. Bhat,
\textit{An index theory for quantum dynamical semigroups}, Trans. Amer. Math. Soc., \textbf{348} (1996) 561--583.

\bibitem{BV05}
J. A. Ball, V. Vinnikov, \textit{Lax-Phillips scattering and
conservative linear systems: a Cuntz-algebra multidimensional
setting}, Mem. Amer. Math. Soc., \textbf{178} (2005).

\bibitem{DG07}
S. Dey, R. Gohm,
\textit {Characteristic functions for ergodic tuples},
Integral Equations and Operator Theory, \textbf{58} (2007), 43--63.

\bibitem{DG11}
S. Dey,; R. Gohm,
\textit {Characteristic functions of liftings},
J. Operator Theory, \textbf{65} (2011), 17--45.

\bibitem{FM78}
E. Fornasini,; G. Marchesini,
\textit {Doubly-indexed Dynamical Systems: State Space Models and Structural Properties},
 Math. Systems Theory, \textbf{12} (1978), 59--72.

\bibitem{GGY08}
J. Gough, R. Gohm, Yanagisawa: \textit{Linear Quantum feedback
Networks}, Phys. Rev. A, \textbf{78} (2008).

 \bibitem{Go04}
R. Gohm, \textit{Noncommutative stationary processes,}
Lecture Notes in Mathematics, 1839, Springer-Verlag, Berlin (2004).

\bibitem{Go10}
R. Gohm,
\textit {Non-commutative Markov chains and multi-analytic operators},
 J. Math. Anal. Appl., \textbf{364} (2010), 275--288.

 \bibitem{Go11}
R. Gohm,
\textit {Transfer function for pairs of wandering subspaces},
Spectral theory, mathematical system theory, evolution equations, differential and difference equations,  385--398, Oper. Theory Adv. Appl., 221,
Birkh\"{a}user/Springer Basel AG, Basel, (2012).


\bibitem{KM00}
B. K\"ummerer, H. Maassen, \textit{A scattering theory for Markov chains},  Infin. Dimens. Anal. Quantum Probab. Relat. Top.  \textbf{3}  (2000), 161--176.


\bibitem{LP67}
P.D. Lax, R.S. Phillips, \textit{Scattering theory}, Pure and Applied Mathematics \textbf{26} Academic press, New York-London, (1967).


 \bibitem{Po89a}
G. Popescu,
\textit{Isometric dilations for infinite sequences of noncommuting operators},
Trans. Amer. Math. Soc., \textbf{316} (1989), 523--536.

\bibitem{Po89b}
G. Popescu, \textit{Characteristic functions for infinite sequences
of noncommuting operators,} J. Operator Theory, \textbf{22} (1989),
51--71.

\bibitem{Po95}
G. Popescu, \textit{Multi-analytic operators on Fock spaces},
Math. Ann., \textbf{303} (1995), 31--46.

\bibitem{Po99}
G. Popescu, \textit{Poisson transforms on some {$C^*$}-algebras
generated by isometries}, J. Funct. Anal.,
 \textbf{161} (1999), 27--61.

\bibitem{Po06}
G. Popescu, \textit{Free holomorphic functions on the unit ball of $\m B(\m H)^n$}, J. Funct. Anal.,
 \textbf{241} (2006), 268--333.

\bibitem{RS79}
M. Reed, B. Simon, Methods of modern mathematical physics. III. Scattering theory. Academic Press [Harcourt Brace Jovanovich, Publishers], New York-London, (1979).

\bibitem{SF70}
B. Sz.-Nagy, C. Foias, \textit{Harmonic analysis of operators on
Hilbert space,} North Holland Publ., Amsterdam-Budapest (1970).

\bibitem{YK03}
M. Yanagisawa, H. Kimura, \textit{Transfer function approach to
quantum control, part I: Dynamics of Quantum feedback systems,}
IEEE Tranactions on Automatic control, \textbf{48} (2003), no. 12,
2107--2120.

\end{thebibliography}
\end{document}